\documentclass[11pt,reqno]{amsart}
\usepackage{amsmath,amsthm,amsfonts,amssymb,amscd}
\usepackage[latin1]{inputenc}

\usepackage{amssymb}
\usepackage{amsthm}
\usepackage{graphicx}
\usepackage{amscd}

\begin{document}

\title{Integrating factors for groups of formal complex
diffeomorphisms}

\author{M. Martelo and B. Scardua}

\subjclass{Primary 37F75, 32S65; Secondary 34M35}

\keywords{Formal diffeomorphism, metabelian group, dicritical map, solvability length}


\date{}
\maketitle

\begin{abstract}
We study groups of formal or germs of  analytic diffeomorphisms in
several complex variables.  Such groups are related to the study of
the transverse structure and dynamics of Holomorphic foliations, via
the notion of holonomy group of a leaf of a  foliation. For
dimension one, there is a well-established dictionary relating
analytic/formal classification of the group, with its algebraic
properties (finiteness, commutativity, solvability, ...). Such
system of equivalences also characterizes the existence of suitable
{\it integrating factors}, i.e., invariant vector fields and one-forms associated
to the group. In this paper we search the basic lines of such
dictionary for the case of several complex variables groups. For
abelian, metabelian, solvable or nilpotent groups we investigate the
existence of suitable formal vector fields and closed differential
forms which exhibit an invariance property under the group action.
Our results are applicable in  the construction of suitable
integrating factors for holomorphic foliations with singularities.
We believe they are a starting point in the study of the connection
between Liouvillian integration and transverse structures of
holomorphic foliations with singularities in the case of arbitrary
codimension.
\end{abstract}

\tableofcontents

\newtheorem{Theorem}{Theorem}[section]
\newtheorem{Corollary}[Theorem]{Corollary}
\newtheorem{Proposition}[Theorem]{Proposition}
\newtheorem{Lemma}[Theorem]{Lemma}
\newtheorem{Sublemma}[Theorem]{Subemma}
\newtheorem{Claim}[Theorem]{Claim}
\newtheorem{Definition}[Theorem]{Definition}
\newtheorem{Example}[Theorem]{Example}
\newtheorem{Question}{Question}
\newtheorem{Problem}[Theorem]{Problem}
\newtheorem{Conjecture}{Conjecture}
\newtheorem{Fact}{Fact}

\newtheorem{Remark}[Theorem]{Remark}
\newcommand{\ep}{\varepsilon}
\newcommand{\eps}[1]{{#1}_{\varepsilon}}

\def\hot{\operatorname{{h.o.t.}}}
\def\GL{\operatorname{{GL}}}
\def\Hom{\operatorname{{Hom}}}
\def\Re{\operatorname{{Re}}}
\def\SL{\operatorname{{SL}}}
\def\Res{\operatorname{{Res}}}
\def\Fol{\operatorname{{Fol}}}
\def\tr{\operatorname{{tr}}}
\def\dim{\operatorname{{dim}}}
\def\Aut{\operatorname{{Aut}}}
\def\GL{\operatorname{{GL}}}
\def\Aff{\operatorname{{Aff}}}
\def\Hol{\operatorname{{Hol}}}
\def\loc{\operatorname{{loc}}}
\def\Ker{\operatorname{{Ker}}}
\def\mI{\operatorname{{Im}}}
\def\Dom{\operatorname{{Dom}}}
\def\Id{\operatorname{{Id}}}
\def\Tni{\operatorname{{Int}}}
\def\supp{\operatorname{{supp}}}
\def\Diff{\operatorname{{Diff}}}
\def\sing{\operatorname{{sing}}}
\def\sing{\operatorname{{sing}}}
\def\codim{\operatorname{{codim}}}
\def\grad{\operatorname{{grad}}}
\def\Ind{\operatorname{{Ind}}}
\def\deg{\operatorname{{deg}}}
\def\rank{\operatorname{{rank}}}
\def\Sep{\operatorname{{Sep}}}
\def\Sat{\operatorname{{Sat}}}
\def\Det{\operatorname{{Det}}}

\section{Introduction and main results}
\label{section:introduction}

The study of groups and germs of complex diffeomorphisms  fixing the
origin is an important tool in Complex Dynamics and in the theory of
Holomorphic Foliations, via the study of holonomy groups (cf.
\cite{Camacho-LinsNeto}) of its leaves. Indeed, the holonomy groups
of (the leaves) of a codimension $n \geq 1$ holomorphic foliation
are (identified with) groups of germs of complex diffeomorphisms
fixing the origin of $\mathbb C^n$. In the codimension $ n =1$ case
these are subgroups of germs of one variable holomorphic maps and
there is a well-established dictionary relating topological and
dynamical properties of (the leaves of) the foliation to algebraic
properties of the group. This is clear in works as \cite{CLNS},
\cite{Ilyashenko}, \cite{Nakai}  and \cite{Scherbakov2}.

All these facts are compiled in some works relating the existence of
suitable ``transverse structures" for the foliation with the
transverse dynamics of the foliation (\cite{Camacho-Scardua},
\cite{Scardua}, \cite{scarduaintegration}). In \cite{Brochero} F.
Brochero studies groups of germs of complex analytic diffeomorphisms
having a fixed point at the origin. For such  groups the author
gives a nice study mainly focused on the analytical or topological
description of the following cases: finite groups, linearizable
groups, abelian groups containing a ``generic" dicritic
diffeomorphism. This important work also motivates some of the
concepts and results in our work.
\subsection{Preliminaries, notation and definitions}
Let us introduce the notation we use throughout this paper.  We denote by $\mathcal O_n$ the ring of
germs at the origin of holomorphic functions of $n$ variables and by
$\hat {\mathcal O_n}$ its
 formal completion\footnote{We refer to the book of D. Eisenbud (\cite{Eisenbud})
 for a detailed construction of the formal completion of the ring $\mathcal
 O_n$,  as well as for some classical  concepts related to the formal objects we deal with in this
 paper.}.

Denote by $\Diff(\mathbb C^n, 0)$ the group of germs of complex
diffeomorphisms fixing the origin $0\in \mathbb C^n$. We denote by
$z=(z_1,...,z_n)$ a system of complex variables in $\mathbb C^n$.
The group of formal diffeomorphisms in $n$ complex variables fixing
the origin is the formal completion $\hat{\Diff}(\mathbb C^n,0)$ of
$\Diff(\mathbb C^n,0)$, obtained from the power series of the
coordinate functions of elements in $\Diff(\mathbb C^n,0)$. This
way, given a formal diffeomorphism $\hat f  \in \hat \Diff(\mathbb
C^n,0)$ we write $\hat f = \hat f^\prime (0) \cdot z +
\sum\limits_{j=2} ^\infty f_j(z)$ where $\hat f^\prime(0) \in \GL(n,
\mathbb C)$ and each $f_j$ is a vector whose coordinates are
homogeneous polynomials of degree $j$ in the variables
$z=(z_1,...,z_n)$. To each formal diffeomorphism $\hat f \in
\Diff(\mathbb C^n, 0)$ we associate its derivative $\hat f
^\prime(0)\in \GL(n,\mathbb C)$. We say that the diffeomorphism
$\hat f$ is {\it tangent to the identity} if $\hat f^\prime
(0)=\Id$. Denote by $\Diff_{\Id}(\mathbb C^n,0)$ the subgroup of
elements tangent to the identity in $\hat{\Diff}(\mathbb C^n, 0)$.
Also put $\Diff_{\Id}(\mathbb C^n,0)=\hat{\Diff}_{\Id}(\mathbb C^n,
0)\cap \Diff(\mathbb C^n,0)$.  This gives  inclusions of
$\Diff(\mathbb C^n, 0) \hookrightarrow\hat \Diff(\mathbb C^n, 0)$
and $\Diff_{\Id}(\mathbb C^n, 0)\hookrightarrow \hat
\Diff_{\Id}(\mathbb C^n, 0)$. A subgroup $G <
\hat\Diff_{\Id}(\mathbb C^n,0)$ is called {\it tangent to the
identity}. The  subgroup of formal diffeomorphisms, {\it tangent to
the identity
 with order $k$}, is defined as ${\hat \Diff}_{k}(\mathbb{C}^{n},0) =
\{\hat f \in {\hat \Diff}(\mathbb{C}^{n},0) \mid \hat f(z) = z +
f_{k+1}(z) + f_{k+2}(z) + \cdots, \, f_{k+1} \ne 0 \}$. Similarly
the group of germs of holomorphic diffeomorphisms at the origin
$0\in \mathbb C^n$, tangent to the identity with order $k$ is
defined as ${\Diff}_{k}(\mathbb{C}^{n},0)={\hat
\Diff}_{k}(\mathbb{C}^{n},0)\cap \Diff(\mathbb C^n,0)$.

The classical theory of groups states the following algebraic
definitions.  Let $G$ be a group. Given elements $\alpha, \beta \in
G$, the {\it commutator} of $\alpha$ and $\beta$ is defined as
$[\alpha, \beta] = \alpha \beta \alpha^{-1} \beta^{-1}$. Given
subgroups $H$ and $L$ of $G$, we define {\it group of commutators}
$[H,L]$ as the subgroup of $G$ generated by the elements of the form
$[\alpha,\beta]$ for $\alpha \in H$ and $\beta \in L$. We put
\[ G^{(0)} = G, \  G^{(n+1)} = [G^{(n)},G^{(n)}]  \ \forall n \geq 0\]
the {\it derived series} of $G$. We define
\[ {\mathcal C}^{0} G = G, \  {\mathcal C}^{n+1} G = [G, {\mathcal C}^{n} G]  \ \forall n \geq 0\]
the {\it descending central series} of $G$. The group $G$ is
solvable if it has a finite derived series. i.e., if $G^{(k)} =  \{
\Id \}$ for some $k \in \mathbb N$. The mininum of such $k$ is called the {\it soluble length} $l(G)$ of $G$. The group $G$ is {\it nilpotent}
if it has a finite descending series, i.e., if
 ${\mathcal C}^{j} G =
\{ \Id \}$ for some $k \in \mathbb N$.

\vglue.1in

\subsubsection{Formal vector fields  and formal diffeomorphisms}

Denote by $\mathcal X(\mathbb C^n,0)$ the $\mathcal O_n$-module of
germs of complex vector fields vanishing at the origin $0\in \mathbb
C^n$ and by $\hat {\mathcal X} (\mathbb C^n,0)$ its formal
counterpart.

\begin{Definition}
\label{Definition:invariantvectorfield}
{\rm  Given a subgroup $G< \hat{\Diff}(\mathbb C^n, 0)$, a formal
vector field $\hat X\in \hat{\mathcal X}(\mathbb C^n,0)$ is
$G$-{\it invariant} if we have $\hat g* \hat X=\hat X, \, \forall \hat g
\in G$. We say that $\hat X$ is {\it projectively invariant} by $G$ if for
each $\hat g \in G$ there is $c_g \in \mathbb C$ such that  $\hat g_* \hat X = c_{\hat g}
\cdot \hat X$.}
\end{Definition}

The Lie algebra $\hat{\mathcal X}_{k}(\mathbb{C}^n,0)$ of formal
vectors fields  of $\mathbb{C}^{n}$ of order $k+1$ is defined by
those vector fields for the form $\hat X=\hat
a_{1}(z)\frac{\partial}{\partial z_{1}} + \cdots + \hat
a_{n}(z)\frac{\partial}{\partial z_{n}}; \hat a_{j}\in \hat{\mathcal
O}_n$, where the minimum order of vanishing of  the  $\hat a_j$ at
the origin is  $k+1$. Thus, $\hat {\mathcal X}_k(\mathbb C^n,0)$ is
the formal completion of ${\mathcal X}_k (\mathbb C^n, 0)$;  the set
of germs of complex analytic vector fields which are singular of
order $k+1$ at $0$. We denote by ${\mathcal X}_{N} (\mathbb C^n, 0)$
the subset of ${\mathcal X} (\mathbb C^n, 0)$ of nilpotent vector
fields, i.e. vector fields whose first jet has the unique eigenvalue
$0$. The formal completion of  $\hat{\mathcal X}_N (\mathbb C^n, 0)$
is denoted by  $\hat{\mathcal X}_{N} (\mathbb C, 0)$.

The expression
\begin{equation}
\label{equ:exp} {\rm exp} (t \hat{X}) = \left({ \sum_{j=0}^{\infty}
\frac{t^{j}}{j!} \hat{X}^{j}(z_{1}), \hdots, \sum_{j=0}^{\infty}
\frac{t^{j}}{j!} \hat{X}^{j}(z_{n}) }\right)
\end{equation}
defines the exponential of $t \hat{X}$ for $\hat{X} \in
\hat{\mathcal X} (\mathbb C^n, 0)$ and $t \in {\mathbb C}$. Let us
remark that $\hat{X}^{j}(g)$ is the result of applying $j$ times the
derivation $\hat{X}$ to the power series $g$. The definition
coincides with the classical one if $\hat{X}$ is a germ of
convergent vector field.

\vglue.1in

We denote by $\Diff_u (\mathbb C^n, 0)$ the subgroup of unipotent
elements of $\Diff(\mathbb C^n, 0)$, more precisely $\varphi \in
\Diff_u (\mathbb C^n, 0)$ if $j^{1} \varphi$ is a unipotent linear
isomorphism (i.e. $j^{1} \varphi - Id$ is nilpotent).
 Analogously we denote $\hat\Diff_u (\mathbb C^n, 0)$
the formal completion of $\Diff_u (\mathbb C^n, 0)$.

According to (\cite{Ilyashenko-Yakovenko}) the exponential map
$\exp:\hat{\mathcal X}_{k}(\mathbb{C}^n,0)\to {\hat
\Diff}_{k}(\mathbb{C}^{n},0)$ is a bijection. Also, it induces a
bijection from $\hat{\mathcal X}_{N} (\mathbb C^n, 0)$ onto
$\hat\Diff_{u}(\mathbb C^n, 0)$.

\begin{Definition}{\rm
Let $\hat f\in \hat{\Diff}_k (\mathbb C^n, 0)$. We denote by $\log
\hat f $ the unique element of $\hat{\mathcal X}_{k} (\mathbb C,0)$
such that $\hat f  = \exp(\log \hat f)$. We say that $\log \hat f$
is the {\it infinitesimal generator} of $\hat f$. Given a map $\hat
h \in \hat{\Diff}_{\Id}(\mathbb C^n,0)$ and $ t \in \mathbb C$ we
denote by $\hat h ^{[t]}=\exp (t.\hat X)$ where $\hat h = \exp \hat
X$.}
\end{Definition}
In general the infinitesimal generator of a (convergent) germ of
diffeomorphism is a divergent vector field (see \cite{Ah-Ro}).

\subsubsection{Solvable length of a Lie algebra}
 Let ${\mathfrak
g}$ be a complex Lie algebra. Given elements $X, Y \in {\mathfrak
g}$ we denote by $[X,Y]$ the Lie bracket of $X,Y \in {\mathfrak g}$.
Given Lie subalgebras ${\mathfrak h}$ and ${\mathfrak l}$ of
${\mathfrak g}$. We define $[{\mathfrak h},{\mathfrak l}]$ the Lie
subalgebra of ${\mathfrak g}$ generated by the elements of the form
$[X,Y]$ for $X \in {\mathfrak h}$ and $Y \in {\mathfrak l}$. We
define
\[ {\mathfrak g}^{(0)} = {\mathfrak g}, \
{\mathfrak g}^{(n+1)} = [{\mathfrak g}^{(n)},{\mathfrak g}^{(n)}]  \
\forall n \geq 0\] the derived series of ${\mathfrak g} $. We define
\[ {\mathcal C}^{0} {\mathfrak g} = {\mathfrak g}, \
{\mathcal C}^{n+1} {\mathfrak g} = [{\mathfrak g}, {\mathcal C}^{n}
{\mathfrak g}]  \ \forall n \geq 0\] the descending central series
of ${\mathfrak g}$.

\begin{Definition}[solvable length of a  Lie algebra, cf. \cite{Martelo-Ribon}]
\label{def:length} {\rm Let  $\mathfrak g$ be a
Lie algebra. We define  $l(\mathfrak g)$ the soluble
length of ${\mathfrak g}$ as

\[ l(\mathfrak g) =
\min \{ k \in {\mathbb N} \cup \{0 \} :  \mathfrak g^{(k)} = \{ Id
\} \}
\]
where $\min \emptyset = \infty$.

We say that $\mathfrak g$ is {\it solvable} if $l(\mathfrak g) <
\infty$.  The Lie algebra  $\mathfrak g$ is {\it nilpotent} if there
exists $j \geq 0$ such that ${\mathcal C}^{j} \mathfrak g = \{ \Id
\}$. If $j$ is the minimum non-negative integer number with such a
property we say that $\mathfrak g$ is of {\it nilpotent class} $j$.}
\end{Definition}

\subsubsection{Lie algebra of a group of diffeomorphisms (cf. \cite{Martelo-Ribon})}

In \cite{Ghys} E. Ghys associates a Lie algebra of formal nilpotent
vector fields to any group of unipotent diffeomorphisms (prop. 4.3
in \cite{Ghys}). In the same spirit we present a construction, from
\cite{Martelo-Ribon},  that associates a non-trivial  Lie subalgebra of
$\hat{\mathcal X} (\mathbb C^n, 0)$ to certain  subgroups $G$ of $\hat
\Diff(\mathbb C^n, 0)$, satisfying some connectedness hypothesis. 
In \cite{Martelo-Ribon} the authors  replace
$G$ with a subgroup $\overline{G}^{(0)}$ of $\hat \Diff(\mathbb C^n,
0)$ containing $G$ that is, roughly speaking, the algebraic closure
of $G$ (with respect to the Krull topology). Such a group satisfies
$l(\overline{G}^{(0)})=l(G)$ (i.e, the  groups $G^{(0)}$ and $G$ have  the same
solvable length) and it has a non-trivial  with analogous algebraic
properties. This construction can be performed to every {\it connected}\,\footnote{A subgroup
$G<\hat{\Diff}(\mathbb C^n,0)$ is
 {\it connected} if the closure in the Zariski
topology of the subgroup of linear parts $j^{1} G < \GL(n,{\mathbb
C})$ of $G$ is connected.} group of formal diffeomorphisms.

\subsubsection*{Associate Lie algebra}
We introduce the concept of {\it associate Lie algebra} of a subgroup of
$\hat{\Diff}(\mathbb C^n,0)$ as presented in \cite{Martelo-Ribon}.
First recall that, according to \cite{Martelo-Ribon} Definition 3.8,
given a  subgroup $G$  of
$\hat{\Diff}(\mathbb C^n,0)$, the {\it
Lie algebra of} $G$ is the  complex Lie subalgebra
${\log (G)}$ of $\hat{\mathcal X} (\mathbb C^n,0)$ given by 
\[
 \log (G)= \{ \hat{X} \in \hat{\mathcal X} (\mathbb C^n,0) : {\rm exp}(t \hat{X}) \in G, \, \
\forall t \in {\mathbb C} \} .
 \]

With this definition the Lie algebra of a discrete group is trivial.
However, under some connectedness hypothesis,  it is possible to  build  a group $\overline{G}^{(0)}\supset G$
with the same algebraic properties of $G$, having a non-trivial Lie algebra $\log(\overline{G}^{(0)})$ and such that,
 $l(\overline{G}^{(0)}) = l(G)$ and   $l(\log(\overline{G}^{(0)})) = l(\log G)$ . Roughly speaking $\overline{G}^{(0)}$ is the algebraic closure of $G$ (with respect to the Krull topology). In particular this construction associates a non-trivial Lie algebra, with analogous algebraic properties, to every  group $G$ of formal diffeomorphisms with connected linear part $DG$ .
In few words, this is done as follows. 
 Let ${\mathfrak m}$ the maximal ideal of
${\mathbb C}[[x_{1},\hdots,x_{n}]]$.
A formal diffeomorphism $\hat h \in \hat{\Diff}(\mathbb C^n,0)$ acts on the 
space   ${\mathfrak m}/{\mathfrak m}^{k+1}$ of $k$-jets in a natural way ($[g] + {\mathfrak m}^{k+1} \mapsto [ g \circ \hat h] + {\mathfrak m}^{k+1}$), defining  an element
$\hat h_{k}$ of $GL({\mathfrak m}/{\mathfrak m}^{k+1})$. 
Let $G$ be a connected subgroup of $\hat{\Diff}(\mathbb C^n,0)$.
Fixed $k \in {\mathbb N}$ we define the group
$C_{k}=\{ \hat h_{k} : \hat h \in G\} \subset GL({\mathfrak m}/{\mathfrak m}^{k+1})$ and the  matrix group $G_{k}$
defined as the smallest algebraic subgroup of $GL({\mathfrak m}/{\mathfrak m}^{k+1})$
containing $C_{k}$.  We define
\[ \overline{G}^{(0)} = \{ \hat h \in \hat{\Diff}(\mathbb C^n,0) : \hat h_{k} \in G_{k} \ \forall k \in {\mathbb N} \} . \]
Then  $\overline{G}^{(0)}$ is a subgroup of $\hat{\Diff}(\mathbb C^n,0)$ containing $G$ that is closed
in the Krull topology. 
Since  a subgroup tangent to the identity is connected we can introduce the following definition: 

\begin{Definition}[Associate Lie algebra, cf. \cite{Martelo-Ribon}]
\label{Definition:asLiealgebragroup} {\rm  Let $G$ a  subgroup of
$\hat{\Diff}_{\Id}(\mathbb C^n,0)$. The {\it
``associate" Lie algebra of} $G$ is defined as  the Lie
algebra of $\overline{G}^{(0)}$. }
\end{Definition}
\begin{Remark}
{\rm It is important to remark that, as observed  in
\cite{Martelo-Ribon}, the Lie algebra of a connected group $G <
\hat{\Diff}(\mathbb C^n,0)$ shares the usual properties of Lie
algebras of Lie groups. This is a consequence of the following
proposition where we can consider $\overline{G}_{0}^{(0)}$ as
the connected component of the identity of $\overline{G}^{(0)}$.

\begin{Proposition}[\cite{Martelo-Ribon} Proposition 3.1]
\label{pro:lie} Let $G < \hat{\Diff}(\mathbb C^n,0)$ be a connected group with associate Lie algebra $\mathfrak{g}$ as in Definition~\ref{Definition:asLiealgebragroup}.
Then ${\mathfrak g}$ is the Lie algebra of $\overline{G}^{(0)}$. The
group $\overline{G}_{0}^{(0)}$ is generated by ${\rm exp}({\mathfrak
g})$. Moreover if $G$ is unipotent then ${\mathfrak g}$ is a Lie
algebra of formal nilpotent vector fields and ${\rm exp}: {\mathfrak
g} \to \overline{G}^{(0)}$ is a bijection.
\end{Proposition}

}
\end{Remark}

We denote by  $\hat{K}_{n}$ the {\it field of fractions} of the ring
of formal power series $\hat {\mathcal O}_n={\mathbb
C}[[z_{1},\hdots,z_{n}]]$.
We have a natural embedding $\hat{K}_n \subset K_n$, where
 by $K_n$ we mean the field of rational functions in $n$ complex variables, i.e., the fraction field of $\mathbb C[z_1,...,z_n]$.
Following the above notions it is then natural to define:

  \begin{Definition}[\cite{Martelo-Ribon} Definition 3.9, dimension of a Lie algebra]
  {\rm 
For a group tangent to the identity $G<
\hat\Diff_{\Id}(\mathbb C^n,0)$, by
   {\it   dimension} of the {\it associate} Lie algebra we mean
 the dimension of the Lie algebra $\log(\overline{G}_{0})$, viewed as
vector space over $\hat{K}_n$.}
\end{Definition}

From now on, by Lie algebra of a group of formal diffeomorphisms tangent to the identity, we shall mean 
its associate Lie algebra.

\subsection{The one-dimensional case}

The construction of a dictionary relating algebraic, dynamical and
analytic  properties of subgroups of diffeomorphisms is a very
important tool in Complex Dynamics and Holomorphic foliations.
  The one-dimensional case has been
addressed by several authors. Below we find a compilation of their
main achievements:

\begin{Theorem}[\cite{Ilyashenko-Yakovenko}]
\label{Theorem:dimensionone}

Let $G < \hat \Diff (\mathbb C,0)$ be a subgroup of one-variable formal complex diffeomorphisms.

\begin{enumerate}

\item $G$ is abelian if, and only if, $G$ is nilpotent if, and only if, $G$ admits a formal invariant
vector field.

\item The following conditions are equivalent:
    \begin{enumerate}
    \item $G$ is
solvable.
 \item $G$ is metabelian, i.e., $G^{(1)}$ is abelian.

\item All elements tangent to the identity in $G$ have a same order of tangency.

\item  $G$ admits a formal vector
field which is projectively invariant by $G$.
\end{enumerate}

\end{enumerate}

\end{Theorem}

 Item (1) above  is essentially a consequence of two
other facts:

\begin{enumerate}
\item[{\rm(i)}] For a subgroup $G < \hat \Diff (\mathbb C,0)$ the {\it
derivative group} $DG=\{\hat g^\prime (0): \, \, \hat g \in G\}$ is
abelian and therefore the  subgroup of commutators $G^{(1)}< G$,
which is  the set of products of commutators in $G$, is {\it tangent
to the identity}, i.e., a subgroup of $G_{\Id}$.

\item[{\rm(ii)}]  Two  elements $\hat f=z + a_{k+1} z^{k+1} + \cdots\in G_{\Id},$ and $ \hat g = z +
b_{\ell +1} z^{\ell +1} + \cdots \in G_{\Id}$, $a_{k+1} . b_{\ell +
1} \ne 0$, tangent to the identity, commute {\it only if} we have
$k=\ell$.

\end{enumerate}

As we shall see, none of the above facts holds for subgroups of
$\hat\Diff(\mathbb C^n,0)$ when $n \geq 2$ (cf. Example~\ref{Example:abelianformalinvariant}). Therefore it is quite
natural to expect that the above mentioned dictionary is much
different or much harder to find, in the $ n \geq 2$ case. To begin
this study is one of the main goals of this work. We also aim on possible
 applications of our results to the framework of holomorphic
foliations (see Section~\ref{section:holonomygroups} and
Proposition~\ref{Proposition:holonomygroups}).

For some of the reasons mentioned above we divide this work in two
parts. The first is mainly, but not only, concerned with the study of  subgroups tangent
to the identity, i.e., groups with all elements tangent to the
identity. The second is about not necessarily  groups tangent to the
identity, but we require the existence of suitable {\it dicritic}
(``radial type") elements in the group.

\subsection{PART I -  Lie algebras of groups and vector fields}
As mentioned above, in the first part we focus on the study of the
Lie algebra of subgroups $G < \hat\Diff (\mathbb C^n,0)$ under the
hypothesis that $G$ is abelian, metabelian or nilpotent.

\subsubsection{Existence of an invariant formal vector field}
In Section~\ref{Section:preliminaries} we prove the following
(compare with (1) in Theorem~\ref{Theorem:dimensionone}):
\begin{Theorem}
\label{Theorem:abelianformalvector} Let  $G< \hat{\Diff}(\mathbb
C^n,0)$ be an abelian subgroup of formal  diffeomorphisms. We have
two possibilities:

\begin{enumerate}
\item $G$ admits an invariant formal vector field.
\item $G_{\Id}=\{\Id\}$, i.e., the only element of $G$ tangent to the identity is the identity.
\end{enumerate}
\end{Theorem}

\begin{Remark}
{\rm With respect to the one-dimensional case we observe:
\begin{itemize}
\item[{\rm(i)}] Unlike the one-dimensional case,
in general the existence of such an invariant vector field is not
enough to assure that the group is abelian (see
Example~\ref{Example:abelianformalinvariant}).

\item[{\rm(ii)}] For $n=1$,   condition (2) above implies that the group is
conjugated to its linear part by a formal diffeomorphism. We believe that this also holds for dimension $n \geq 2$.
\end{itemize}}
\end{Remark}

\subsubsection{Abelian subgroups and nilpotent subalgebras}

 In \cite{Brochero} it is proved
(Proposition 4.1) that {\it every nilpotent subalgebra $\mathcal L$
of $\hat{\mathcal  X}(\mathbb C^2, 0)$ is metabelian}. As a
consequence, if $\mathcal R$ is the center of $\mathcal L$ then
$\mathcal R \otimes {\hat K}_2$ is a vector space of dimension $1$
or $2$ over $\hat { K}_2$. According to  \cite{Brochero}
Corollary~4.4, if the dimension is one then there is a formal vector
field $\hat X \in \mathcal X(\mathbb C^2,0)$ such that $\exp(\hat
X)\in G$ and such that for each element $g \in G$ there exists a
rational (meromorphic) function $N_g\in K_2$ with $\hat X(N_g)=0$ and
$g=\exp(N_g \hat X)$. In case the dimension is $2$ there are
(commuting) maps $\hat f,  \hat g \in G$ such that $G< \{ \hat f
^{[t]} \circ \hat g ^{[s]}, \, s, t \in \mathbb C\}$. Using this we
obtain the following rephrasing of Brochero's result:

\begin{Theorem}[cf. \cite{Brochero}]
\label{Theorem:A} Let $G < {\Diff}_{\Id}(\mathbb{C}^{2},0)$ be
an abelian subgroup tangent to the identity. We have the following
possibilities:
\begin{enumerate}

\item[{\rm(i)}] $ G$ leaves invariant an exact rational one-form, say
$\hat \omega=dT$ for some rational function $T \in K_2$.

\item[{\rm(ii)}] $ G$ embeds into the flow of a  formal  vector field $\hat {X} \in \log(G)$.

\item[{\rm(iii)}]  There are two
invariant linearly independent commuting  formal vector fields $\hat X, \hat Y \in \log(G)$.
\end{enumerate}
\end{Theorem}

Though this statement is essentially already contained in \cite{Brochero}, we give an alternative proof which will indeed allow an extension of this result for higher dimension $n \geq 2$ and we shall prove in Section~\ref{Section:nilpotent}:

\begin{Theorem}
\label{Theorem:nilpotent} Every nilpotent subalgebra $\mathfrak{l}$
of $\hat{\mathcal{X}}(\mathbb{C}^n,0)$  has length $l(\mathfrak{l})$
at most $n$. As a consequence, if $G < \hat{\Diff}(\mathbb C^n,0)$
is a nilpotent group then the solubility length of $G$ satisfies
$l(G) \leq n$.
\end{Theorem}

As a consequence we obtain, the following  immediate corollary of (the proof of) Theorem~\ref{Theorem:nilpotent}:
\begin{Corollary}
\label{Corollary:abeliandimensionn} Let $G <
\hat{\Diff}_1(\mathbb{C}^{n},0)$ be abelian subgroup tangent to the
identity, with Lie algebra $\mathfrak{g}$. Then the possibilities are:
\begin{enumerate}

\item There are $\hat X_1,...,\hat X_n\in \mathfrak{g}$ such that $[\hat {X}_j,\hat{X}_r] = 0, \forall i,j$ and
$G < \langle \exp(t_1\hat X_1)\circ\cdots\circ
\exp(t_n \hat {X}_n)\mid t_j\in\mathbb{C}\rangle$.

\item There are $\hat X_1,...,\hat X_l\in \mathfrak{g}$, for some $l\in\{1,\ldots,n - 1\}$,
    such that $[\hat {X}_j,\hat{X}_r] = 0, \forall i,j$ and
\[
G < \langle  \exp(u_1\hat X_1)\circ\cdots\circ
\exp(u_l\hat {X}_l)\mid u_j \in K_n,\, \,  \hat {X}_r(u_j) = 0,  \,\forall \, 1\leq
r,j\leq l\,\rangle.
\]

\end{enumerate}
\end{Corollary}

\subsubsection{Abelian groups and closed one-forms}

Before stating our next main results we observe
that in some main applications of the results in
Theorem~\ref{Theorem:dimensionone} (case $n=1$), the winning
strategy is to construct from the information on the holonomy groups
of the foliation, some suitable differential forms which allow to
``integrate" the foliation (as for instance a foliation admitting a
Liouvillian first integral) (see \cite{Camacho-Scardua} for instance). More precisely, in dimension $n=1$ a
formal vector field $\hat {X} \in \hat{\mathcal X} (\mathbb C,0)$
can be written either as a linear vector field or, in case it has a
zero of order $\geq 2$ at the origin, as
\[
\hat {X}(z) = \frac{z^{k+1}}{ 1 + \lambda z^k} \frac{d}{dz}
\]
for some $\lambda \in \mathbb C$ and $k \in \mathbb N$. Let us focus
on the non-linear case. The {\it duality equation} $\hat \omega
\cdot \hat {X}=1$ has,  in this dimension one case, a single
solution
\[
\hat \omega= \lambda\frac{dz}{z} + \frac{dz}{z^{k+1}}
\]
This expression, is the expression of general closed meromorphic
one-form with an isolated pole of order $k+1$ at the origin
$0\in\mathbb C$, residue $\lambda$, in  a suitable coordinate
system. It is a particular case of the so called Integration Lemma
(see for instance \cite{Scardua} Example 1.6 page 174, or
Proposition~\ref{Proposition:integrationlemma} in
Section~\ref{section:integrationlemma}).

\underline{In dimension $n=1$}, given   a formal diffeomorphism
$\hat g \in  \hat \Diff(\mathbb C,0)$ and $\hat {X}$ and
 $\hat\omega$ satisfying the duality equation as above, we have:
\begin{enumerate}
\item $\hat g _* \hat {X} = \hat {X}\iff \hat g^* \hat \omega = \hat\omega$

\item $\hat g_* \hat {X} = c_{\hat g} \hat {X}$ for some $c_{\hat g} \in \mathbb C^*
\iff g^* \hat \omega = \frac{1}{c_{\hat g}}
\hat \omega$
\end{enumerate}

Finally, notice that, in dimension-one each formal or meromorphic
one-form is closed. This suggests, in view of
Theorem~\ref{Theorem:dimensionone}  and all the above, that one may
expect to obtain results relating algebraic properties of subgroups
of $\hat \Diff (\mathbb C^n,0)$ with the existence of suitable
{\em closed} one-forms.

 \vglue.1in We shall see (cf. Theorem~\ref{Theorem:B} below) that {\it a subgroup $G
< \hat{\Diff}(\mathbb C^2,0)$ admitting two commuting formal
invariant vector fields, exhibits two independent invariant  closed formal
meromorphic one-forms}. A further investigation of this situation involves the following notion:

\begin{Definition}[Formal separatrices]
\label{Definition:formalcurve} {\rm A {\it formal curve} of $2$
complex variables is defined as follows:  In the local ring
$\hat{\mathcal O_2}$ we introduce the equivalence relation $\hat f
\thicksim\hat g\iff \hat \varphi=\hat u. \hat \psi$ for some unit
$\hat u\in \hat{\mathcal O_2}$, i.e., for some power series $\hat u$
with first coefficient $u_0\ne 0$. By a {\it formal curve} we mean
an equivalence class of a function $\hat \varphi\in \hat{\mathcal
O_2}$ that satisfies $\hat \varphi(0)=0$, that is, a non-invertible
formal power series. Such a formal curve is called {\it invariant}
by a formal complex diffeomorphism $\hat f\in \hat \Diff(\mathbb
C^2,0)$ if ${\hat f}^*\varphi=\varphi \circ \hat f$ is equivalent to
$\varphi$ in the above sense. Such a formal curve will be called a
{\it separatrix} of a subgroup $ G < \hat \Diff(\mathbb C^2, 0)$ if
it is invariant by each element  of this group. The {\it tangent
space} of a formal curve with representative $\hat \varphi$ is
defined as the linear subspace of $\mathbb C^2$ given by the kernel
of $D \hat \varphi(0)\colon \mathbb C^2 \to \mathbb C$. Two formal
curves with representatives $\hat \varphi$ and $\hat \psi$ are
called {\it transverse} if:
\begin{enumerate}

\item Each tangent space has dimension one.

\item The tangent spaces span $\mathbb C^2$.
\end{enumerate}}

\end{Definition}

As a converse of (iii) in Theorem~\ref{Theorem:A}
 we have:
\vglue.1in \begin{Theorem} \label{Theorem:B} Let $ G < {\hat
\Diff}(\mathbb{C}^2,0)$ be a subgroup admitting two linearly independent  invariant
commuting  formal vector fields.
Then  $ G$ admits two closed,
independent, formal meromorphic, invariant one-forms.
The group $G$ is abelian provided that one of the following conditions is satisfied:

\begin{enumerate}
\item
 $ G$  is
tangent to the identity.

 \item $G$  exhibits two formal transverse
separatrices.
\end{enumerate}
\end{Theorem}
 The notions of formal closed meromorphic one-form and
other formal objects are clearly stated in
Section~\ref{section:flatabelianclassification}
(Definition~\ref{Definition:formalobjects}). As a spolium  of the
proof of the second part of Theorem~\ref{Theorem:B} we obtain normal
forms for abelian subgroups admitting two transverse separatrices
and having Lie algebra of dimension two (cf.
Remark~\ref{Remark:normalformdouble}).

\subsubsection{Metabelian groups}
All the above is concerned with the abelian case. As for the
metabelian  non-abelian case we have:

\begin{Theorem}
\label{Theorem:metabelian}
 Let $G <
\hat{\Diff}(\mathbb{C}^{2},0)$ be a metabelian non-abelian subgroup.
Assume that the group of commutators   $G^{(1)}=[G,G]$ is tangent to
the identity, for instance if the derivative group $DG<
\GL(2,\mathbb C)$ is abelian. Denote by $l( G^{(1)})$
  the Lie algebra of  $ G^{(1)}=[ G, G]$. We have the
following possibilities:
\begin{enumerate}
\item[{\rm(i)}]
$l( G^{(1)})$ is one-dimensional. There  is a formal vector field
$\hat {X}$  with $\exp \hat X \in G$,  such that, for each  $\hat
g\in G$ there is a rational function $T_{\hat g}\in K_2$ that satisfies:
$\hat {X} (T_{\hat g})=0$ and $\hat g^*(\hat {X}) = T_{\hat g} \cdot
\hat {X}$.
\item[{\rm(ii)}]
$l( G^{(1)})$ is two-dimensional. There
 are two $\mathbb C$-linearly independent formal vector fields $\hat {X},
 \hat {Y} \in \mathfrak g=\log G$,
  such that
 \begin{itemize}
 \item[{\rm(ii.1)}] $[\hat {X}, \hat {Y}]=0$

\item[{\rm(ii.2)}] For each  $\hat g\in G$ there are
$(s_1,t_1),(s_2,t_2)\in\mathbb{C}^2$
 $\mathbb C$-linearly independent that satisfy $\hat g^*(\hat {X}) =
 s_1\hat {X} + t_1\hat {Y}$ and $\hat g^*(\hat {Y}) = s_2\hat {X} + t_2\hat {Y}$.
 \end{itemize}
 \end{enumerate}
 Furthermore, in this last case  there are two $\mathbb C$-linearly independent
closed formal meromorphic one-forms $\hat \omega_j$, $(j = 1,2)$ and
$a_j, b_j\in\mathbb{C}^*$ such that
\[
 \hat g^*(\hat \omega_j) = a_j\hat \omega_1 +
b_j\hat \omega_2, \, \,\forall \hat g \in  G.
\]

\end{Theorem}

Groups as in (ii.2) above are studied in
Section~\ref{section:metabelian} (cf.
Remark~\ref{Remark:normalformdouble}).

\subsection{PART II - Groups containing dicritic  diffeomorphisms}

The second part of this work is dedicated to the study of subgroups
of formal diffeomorphisms under the hypothesis of existence of a
suitable dicritic (radial type) element. More precisely, according
to \cite{Brochero},   a diffeomorphism  $\hat f\in
\hat{\Diff}_{k}(\mathbb{C}^{n},0)$ is called {\it dicritic} if $
\hat f(z) = z + \hat f_{k + 1}(z) + \hat f_{k + 2}(z) +\cdots$,
where $\hat f_{k + 1}(z) = f(z)z$ and $f$ is a homogeneous
polynomial of degree $k$. A formal vector field $\hat {X}\in
\hat{\mathcal X}_{k}(\mathbb{C}^n,0), \, k\geq 1$ is called {\it
dicritic} if $\hat {X} = f(z)\vec{R} + (p_{k + 2}^{(1)} +
\cdots)\frac{\partial}{\partial z_{1}} + \cdots + (p_{k + 2}^{(n)} +
\cdots)\frac{\partial}{\partial z_{n}}$ where $f_\nu$ is a
homogeneous polynomial of degree $\nu$ and $f$ is homogeneous of
degree $k$ and $\vec{R} = z_{1}\frac{\partial}{\partial z_{1}} +
\cdots + z_{n}\frac{\partial}{\partial z_{n}}$ is the radial vector
field.

\begin{Definition}[Regular dicritic vector fields and diffeomorphisms]
{\rm  We also introduce we the following  useful subclasses of
dicritic diffeomorphisms and vector fields. A formal vector field
$\hat{\mathcal X}_{k}\ni \hat {X}= f(z)\vec{R} + (p_{k + 2}^{(1)} +
\cdots)\frac{\partial}{\partial z_{1}} + \cdots + (p_{k + 2}^{(n)} +
\cdots)\frac{\partial}{\partial z_{n}}$ is called {\it regular
dicritic} if $\hat {X}$ is dicritic and there are $i_{0}, j_{0}\in
\{1,\ldots,n\}$, such that $f$ and $z_{j_{0}}p_{k + 2}^{(i_{0})} -
z_{i_{0}}p_{k + 2}^{(j_{0})}$ are coprime. {\it This implies that
$0$ is an isolated singularity\footnote{Indeed, this is equivalent
to the fact that the singularity is isolated, after a number of
blowing-ups at the vector field.} of $\hat {X}$.} A formal
diffeomorphism  $\hat f\in \hat{\Diff}_k(\mathbb C^n,0)$ is called
{\it regular dicritic} if its infinitesimal generator is a regular
dicritic vector field. }
\end{Definition}

We shall say that a subgroup $G< \hat \Diff(\mathbb C^n,0)$ is {\it quasi-abelian} if
 its subgroup $G_{\Id}$ of  elements tangent to the identity is abelian.

Our next results are analogous to those in
Theorem~\ref{Theorem:dimensionone}, for groups containing a regular
dicritic  element. \vglue.1in
\begin{Theorem}
\label{Theorem:C}   For a  subgroup $G< \hat{\Diff}(\mathbb C^n, 0)$ containing
a regular dicritic diffeomorphism, the following conditions are equivalent:

 \begin{enumerate}
 \item $G$ is quasi-abelian.
 \item $G$ admits a projectively invariant regular dicritic formal vector
field.
\end{enumerate}
\end{Theorem}

In particular we obtain:

\begin{Corollary}
\label{Corollary:dicriticabelianvector} A  subgroup of $\hat
\Diff(\mathbb C^n,0)$ tangent to the identity and containing a
regular dicritic diffeomorphism is abelian if and only if it admits
an invariant formal vector field.
\end{Corollary}

The proof of  Theorem~\ref{Theorem:C}  also shows that
\begin{Proposition}
\label{Proposition:exttheorem2}
Let $ G < {\hat
\Diff}(\mathbb{C}^{n},0)$ be a subgroup containing a regular dicritic
diffeomorphism $\hat f\in G $ and containing its derivative subgroup,  $D G < G $.
The following conditions are equivalent:
\begin{enumerate}

\item $ G $ is abelian

\item  $D G $ is abelian and $\hat X = \log \hat f$ is
invariant by $G$.
\end{enumerate}

\end{Proposition}

 Regarding the case of
metabelian groups  we have:

\begin{Theorem}
\label{Theorem:metabeliandicritic} A subgroup of formal
diffeomorphisms containing a regular dicritic  diffeomorphism and
with abelian derivative group is metabelian provided that it admits
a projectively invariant formal vector field.
\end{Theorem}

As a  partial converse of Theorem~\ref{Theorem:metabeliandicritic}
we have:

\begin{Proposition}
\label{Proposition:converseTheoremB} Let $ G < {\hat
\Diff}(\mathbb{C}^{n},0)$ be a metabelian subgroup containing a
regular dicritic diffeomorphism with order of tangency $k$. Suppose
that $D G $ is abelian and there is a linear diffeomorphism $\hat h
\in G $, given by $\hat h(z) = \lambda \cdot \Id$, where $\lambda\in
\mathbb C$  is such that  $\lambda^{k} \neq 1$, $\lambda^{k + 1}
\neq 1$. Then there is a formal vector field $\hat {X}\in
\hat{\mathcal X}_{j}(\mathbb{C}^{n},0)$, $j\geq 2$, which is
projectively invariant by $G$.
\end{Proposition}

 As an  application we study the case where a group with two generators,
one of which is linear, is metabelian (cf.
Corollary~\ref{corollary:theoremB}).

 Next we state an
equivalence similar to the dimension one case, but for groups that
contain some dicritic diffeomorphism. Theorem~\ref{Theorem:D} below
is related to Theorem 4.1 and Corollary 4.2 in \cite{Brochero} and
to our Example~\ref{Example:solvnotmeta} of a solvable  group of
formal diffeomorphisms  tangent to the identity  which is not
metabelian. This example shows the need of our assumption of
existence of a dicritic element in the subgroup of commutators  in any
extension of Theorem~\ref{Theorem:dimensionone} to higher dimension.
\vglue.1in
\begin{Theorem}
\label{Theorem:D}
 For a  subgroup $G< \hat{\Diff}_{\Id}(\mathbb C^n,0)$ of formal diffeomorphisms tangent to the identity
 and containing a dicritic diffeomorphism tangent to the identity with  order  $k$,
the following statements are equivalent:
\begin{enumerate}
\item[{\rm(1)}]
The group is abelian.
\item[{\rm(2)}]
The group  is nilpotent.
\item[{\rm(3)}]  Every nontrivial element in the group
is tangent to the identity with  order $k$.
\end{enumerate}
\end{Theorem}

\vglue.1in

Our results apply to the study of foliations on complex projective
spaces and other ambient manifolds as well. The class of
singularities which correspond, via the holonomy of its
separatrices, to the class of regular dicritic diffeomorphisms is to
be formally introduced and studied in a forthcoming work. Using an
adaptation of a classical result due to Hironaka and Matsumara
(\cite{[Hir-Mat]}, \cite{[Hart]}) we may be able to move from the
``formal world" (considered is this paper) to the
``analytic/convergent world", which is the natural ambient to the
study of holomorphic foliations with singularities.

A final word should be said about the possible applications of our
results. We are interested in the study of Liouvillian integration
for holomorphic foliations of codimension $n \geq 1$. As suggested
by the codimension one cases (see for instance \cite{ScarduaJPAA}),
this passes through the comprehension of algebraic, geometric and
formal structures of subgroups of $\hat \Diff(\mathbb C^n,0)$ in
terms  we propose in this work.

\vglue.3in

 \centerline{ \textsc{Part I - Groups of diffeomorphisms
tangent to the identity}}
\section{Lie algebras of groups and vector fields}
\label{Section:preliminaries}

Now  we proceed to prove Theorem~\ref{Theorem:abelianformalvector}.
   Some steps in the proof of
the following  well-known  lemma will be used later on this paper:

\begin{Lemma}
\label{fluxo} If $\hat f\in\hat \Diff(\mathbb C^n,0)$ commutes with
the time one flow map  of a formal vector field $\hat {X}\in
\hat{\mathcal X}_{k}(\mathbb C^n,0)$, $k\geq 2$ then $\hat f$
commutes with the flow of $\hat {X}$ for all time $t\in\mathbb{C}$.
\end{Lemma}

\begin{proof}
Let $\hat \Phi _{t}$ be the (formal) flow of $\hat {X}$, which is
defined by $\hat \Phi_t:=\exp(t\hat {X})\in \hat\Diff(\mathbb
C^n,0)$.  Then $\hat \Phi _{1}\circ \hat f = \hat f \circ \hat \Phi
_{1}$. We claim that $\hat \Phi _{t}\circ \hat f = \hat f \circ \hat
\Phi _{t}$ for all $t\in \mathbb{Z}$. First we prove this by
induction, $\forall t\in \mathbb{N}$. In fact, this is true for $t =
1$, Suppose that equality holds for $n\in\mathbb{N}$. Then
\[
\hat \Phi _{n + 1}\circ \hat f = \hat \Phi \circ \hat \Phi _{n}\circ
\hat f  = \hat \Phi\circ \hat f\circ\hat \Phi _{n} = \hat f \circ
\hat \Phi\circ\hat \Phi _{n} = \hat f \circ \hat \Phi _{n + 1}
\]
thus $\hat \Phi _{t}\circ \hat f = \hat f \circ \hat \Phi _{t}$,
$\forall t\in\mathbb{N}$. On the other hand, if $\hat \Phi _{t}\circ
\hat f = \hat f \circ \hat \Phi _{t}$ then $\hat \Phi _{-t}\circ
\hat f = \hat f \circ \hat \Phi _{-t}$. Consequently $\hat \Phi
_{t}\circ \hat f = \hat f \circ \hat \Phi _{t}$, $\forall
t\in\mathbb{Z}$. Now to show that $\hat \Phi _{t}\circ \hat f = \hat
f \circ \hat \Phi _{t}$, $\forall t\in\mathbb{C}$, is sufficient to
prove this equality in the spaces of jets, i.e. in
$\mathcal{J}^k(\mathbb{C}^n,0) = \mathbb{C}[[z]]/\mathfrak{m}^{k +
1}$ $($this has a natural identification with the space of
polynomials of degree less than or equal to $k )$, where
$\mathfrak{m} = \{ \hat f\in \mathbb{C}[[z]]/ \hat f(0) = 0\}$ is
the maximal ideal of $\mathbb{C}[[z]]$. Indeed, given $k\in
\mathbb{N}$ we have that $j^k\circ\hat \Phi _{t}\circ \hat f =
(f_{1},\ldots,f_{n})$, where the truncation of formal series
$$j^k:\mathbb{C}[[z]]\rightarrow \mathcal{J}^k(\mathbb{C}^n,0),$$ is
defined by $j^k(\hat f) = \hat f \mod \mathfrak{m}^{k + 1}$, we have
that $$f_{l} (z) = \sum\limits_{|N|\leq k} P_{N}^{l}(t)z^{N}$$ e
$P_{N}^{l}(t)$ is a polynomial of degree less than or equal to
$|N|$. Similarly, we have\\ $j^k\circ \hat f\circ\hat \Phi _{t} =
(\widetilde{f}_{1},\ldots,\widetilde{f}_{n})$ where
$$\widetilde{f}_{l} (z) = \sum \widetilde{P}_{N}^{l}(t)z^{N}$$ e
$\widetilde{P}_{N}^{l}(t)$ is a polynomial of degree less than or
equal to $|N|$. now, as $\hat \Phi_{t}\circ \hat f = \hat f \circ
\hat \Phi _{t}$, $\forall t\in\mathbb{Z}$, for each
$N\in\mathbb{N}^n$ with $|N|\leq k$, we have that
$P_{N}^{l}(t)\mid_{\mathbb{Z}} =
\widetilde{P}_{N}^{l}(t)\mid_{\mathbb{Z}}$, now as these are
polynomial and coincide in $\mathbb{Z}$, we have that $$P_{N}^{l}(t)
= \widetilde{P}_{N}^{l}(t), \forall t\in \mathbb{C}.$$ in
consequence $f_{l}(z) = \widetilde{f}_{l}(z)$  $\forall X\in\mathbb
C^n$, $l\in\{1,\ldots,n\}$. therefore $j^k\circ\hat \Phi _{t}\circ
\hat f = j^k\circ \hat f\circ\hat \Phi _{t}$, $\forall
t\in\mathbb{C}$ e $k\in\mathbb{N}$. So $$\hat \Phi _{t}\circ \hat f
= \hat f \circ \hat \Phi _{t}, \forall t\in\mathbb{C}.$$
\end{proof}

\begin{proof}[Proof of Theorem~\ref{Theorem:abelianformalvector}]
We may assume that $ G_{\Id}$ is nontrivial. Thus there is $\hat
f\in G_{1}$ which is of the form $\hat f=\exp(\hat {X})$ for some
formal vector field $\hat {X}\in \hat{\mathcal X}_{j}(\mathbb
C^n,0)$, $j\geq 2$. Since $ G$ is abelian, for any $\hat g\in G$,
$\hat g \circ  \hat f(z) = \hat f \circ \hat g(z)$, {\it i.e}, $\hat
g \circ \exp(\hat {X})(z) = \exp(\hat {X})\circ \hat g(z)$. Thus,
from the previous lemma, we have $\hat g \circ  \exp(t\hat {X})(z) =
\exp(t\hat {X})\circ \hat g(z)$, $\forall t\in\mathbb{C}$ or
equivalently $\hat g \circ \exp(t\hat {X})\circ \hat g^{-1} =
\exp(t\hat {X})$, $\forall t\in\mathbb{C}$. Therefore $\hat g^*\hat {X} = \hat {X}$, $\forall \hat g\in G$.
 \end{proof}

 \begin{Remark}
 \label{Remark:complementtoTheoremabelian}
 {\rm Regarding case (2) in Theorem~\ref{Theorem:abelianformalvector}
 we observe that, in $DG$ is abelian and algebraic (see definition in \cite{Rosenlitch} or Remark~\ref{Remark:linearalgebraic}) and
the identity is the only  element tangent to the identity,  the map
$\hat g \mapsto D\hat g(0)$ gives an (abstract)  group isomorphism $
G\cong D G$. Now, according to Remark~\ref{Remark:linearalgebraic}
either $ DG$ is finite (and therefore
analytically conjugated to a finite group of diagonal periodic
linear maps) or the Zariski closure ${\overline{ DG}}$ contains a
linear flow. In this last case, as remarked above, there is a (linear) vector
field $\hat {X}$ which is invariant under the action of  $DG$.
}
\end{Remark}

 \begin{Example}
 \label{Example:abelianformalinvariant}
 {\rm
Now we give some examples showing that the conditions in our main
results, cannot be dropped.
 \begin{enumerate}

 \item The converse of Theorem~\ref{Theorem:abelianformalvector} is
 not always true
  for dimension ($n \geq 2$). In fact if $\hat f(x,y) = (2x,4y)$
  and $\hat g(x,y) = (x,x + y)$ then $ G = \langle \hat f, \hat g\rangle$
  is not abelian, however, $ G$ is invariant by $\hat {X}$,
   where $\exp(\hat {X}) = (x,y + x^2)$. As for the tangent to the identity case, let  $\hat f(x,y) =
   \exp(x^2y\frac{\partial}{\partial x})$
  and $\hat g(x,y) = \exp(x^3y^2\frac{\partial}{\partial x})$, then $G=<\hat f, \hat g>$
  is not abelian, however, $G$ is invariant by
  $\hat{X} = -xy\frac{\partial}{\partial x} + y^2\frac{\partial}{\partial y}$.

   \item In dimension $k = 1$,
   we have that a group $G < {\hat \Diff}_1(\mathbb{C},0)$ of
diffeomorphisms tangent to the identity is abelian if and only if
there is a formal vector field $\hat {X}\in \hat{\mathcal
X}_{k}(\mathbb{C},0)$ $(k\geq 1)$, such that $G < \langle \exp(t\hat
{X})\mid t\in\mathbb{C}\rangle$.  For $(n \geq 2)$ if there is a
formal vector field $\hat {X}\in \hat{\mathcal X}_{k}(\mathbb
C^n,0)$ $(n\geq 2)$, such that $G < \langle \exp(t\hat {X})\mid
t\in\mathbb{C}\rangle$ then $ G$ is abelian, however again  the
converse is not always true. This is due to the fact that for $n =
1$, if the Lie bracket of two vector field $\hat {X}\in
\hat{\mathcal X}_{k}(\mathbb{C},0)$ and $\hat {Y}\in \hat{\mathcal
X}_{r}(\mathbb{C},0)$ is zero $([\hat {X},\hat {Y}] = 0)$ then $r =
k$ and there is $c\in\mathbb{C}^{*}$ such that $\hat {X} = c\, \hat
{Y}$. However this last fact is not always true in dimension $n \geq
2$  as can be seen in the following examples:
\begin{enumerate}
\item Let $a\in\mathbb{C}^{*}$ be constant and take $\hat {X}(x,y)
= xy\frac{\partial}{\partial x} - y^{2}\frac{\partial}{\partial y}
\, , \, \, \hat {Y}(x,y) = ax^{2}y^{2}\frac{\partial}{\partial x} -
axy^{3}\frac{\partial}{\partial y}.$
\item $\hat {X}(x,y) = (x^{2} + 3xy)\frac{\partial}{\partial x} +
(3xy + y^{2})\frac{\partial}{\partial y}, \, \, \hat {Y}(x,y) =
 (3x^{3} - 5x^{2}y + xy^{2} + y^{3})\frac{\partial}{\partial x} +
(x^{3} + x^{2}y - 2xy^{2} + 3y^{3})\frac{\partial}{\partial y}.$
\item

For $k\geq 1$, we have: $\hat {X} = (x^{k +
1})\frac{\partial}{\partial x} + (x^k.y)\frac{\partial}{\partial y}
\, , \, \, \hat {Y} = (y^k.x)\frac{\partial}{\partial x} + (y^{k +
1})\frac{\partial}{\partial y}$

\item $\hat {X}(x,y) = x^{2}\frac{\partial}{\partial x} +
xy\frac{\partial}{\partial y}, \, \,
\hat {Y}(x,y) = xy\frac{\partial}{\partial x} +
y^{2}\frac{\partial}{\partial y}.$
\end{enumerate}
\end{enumerate}

} \end{Example}

\section{Holonomy groups and algebraic abelian groups} In the treatment of abelian
groups we must pay some attention to the linear case. Given a
subgroup $G< \hat{\Diff}(\mathbb C^n,0)$, the derivative map
$D\colon \hat \Diff(\mathbb C^n,0)\to \GL(\mathbb C,n), \, \hat f
\mapsto D\hat f:={\hat f}^\prime(0)$,  induces  by restriction to
any  subgroup $ G < \hat \Diff(\mathbb C^n,0)$ a homomorphism
$D\colon  G \to \GL(\mathbb C,n)$. The kernel of this homomorphism
 is the subgroup $ G _{\Id}$  and the image is the
 {\it derivative subgroup} $D G < \GL(\mathbb C,n)$. If $ G_{\Id}$
 is trivial then we have an injective group homomorphism $G \hookrightarrow \GL(\mathbb C,n)$.

\begin{Remark}[Linear algebraic groups]
\label{Remark:linearalgebraic} {\rm A complex {\it linear algebraic}
group is a subgroup of the group of invertible $n\times n$ complex
matrices (under matrix multiplication) that is defined by complex
polynomial equations. Let therefore $G< \GL(n,\mathbb C)$ be an
infinite linear algebraic group. Then its Lie Algebra $\mathcal
L(G)$ is not trivial and we may choose  a (linear)  vector field
$X\in l(G)$. The Zariski closure $\overline{ \{X_t\}_{t \in \mathbb
C} }$ of the flow $\{X_t\}_{t \in \mathbb C} < \GL(n,\mathbb C)$ of
$X$ in $\mathbb C^n$ is a closed abelian  subgroup of the closure
$\bar G$. Since $\overline {\{X_t\}_{t\in \mathbb C} }$ is abelian,
there is a closed one-parameter subgroup $H$, it contains a
one-dimensional linear algebraic subgroup of $G$.}
\end{Remark}

  Let now $G < \hat{\Diff}(\mathbb{C}^{2},0)$ be an
abelian subgroup. The Lie algebra of $ G_{\Id}=  G \cap
\hat{\Diff}_{\Id}(\mathbb{C}^{2},0)$ has dimension $\leq 2$. If the
dimension is zero then $ G_{\Id}=\{\Id\}$ and the map $ G \to
\GL(2,\mathbb C)$ embeds $ G$ into an abelian linear group. If
moreover the group $G$ is  algebraic, then the derivative group $DG<
\GL(n,\mathbb C)$ is also algebraic and,  in the case
$G_{\Id}=\{\Id\}$, we have from Remark~\ref{Remark:linearalgebraic}
that, either $ G$ is finite or its image in $\GL(2,\mathbb C)$
contains a linear flow in its closure.

From Remark~\ref{Remark:linearalgebraic} and  Theorem~\ref{Theorem:C} we have:

\begin{Corollary}
\label{Corollary:commutative} Let $G <
\hat{\Diff}(\mathbb{C}^{2},0)$ be an algebraic  commutative
subgroup\footnote{The group $G$ is not necessarily tangent to the identity.}. The possibilities  are:
\begin{itemize}

\item[{\rm(i)}] The Lie algebra of $G_{\Id}$ has  dimension zero and the possible cases are:

     \begin{enumerate}
     \item $ G$ is finite and therefore analytically linearizable.

\item $ G$  embeds into a linear flow.

\end{enumerate}
\item[{\rm(ii)}] The Lie algebra of $G_{\Id}$ has  dimension one:
$ G_{\Id}$ leaves invariant an exact rational one-form, say $\hat
\omega=dT$ for some rational function $T\in K_2$.

\item[{\rm(iii)}] The Lie algebra of $G_{\Id}$ has dimension two:
$ G_{\Id}$ admits two closed, independent, formal meromorphic,
invariant one-forms.
\end{itemize}

\end{Corollary}

Next remark we show how
algebraic groups may appear in our framework.

\subsection{Holonomy groups of holomorphic foliations}
\label{section:holonomygroups}
 A dimension one {\it holomorphic
foliation} with singularities ${\mathcal F}$ of a complex manifold
$M$ is defined by a pair $({\mathcal F}^\prime, S)$ where $S\subset
M$ is a codimension $n$ analytic subset of $M$ and ${\mathcal
F}^\prime$ is a dimension one holomorphic foliation of the open
manifold $M^\prime=M \setminus S$, in the usual sense. We may choose
the subset $S\subset M$ as minimal in the sense that ${\mathcal
F}^\prime$ admits no extension as a (nonsingular) foliation to an
open subset intersecting $S$ and in this case we say that
$S=\sing({\mathcal F})$ is the singular set of ${\mathcal F}$. This
is a discrete subset of $M$. A {\it leaf} of ${\mathcal F}$ is by
definition a leaf of ${\mathcal F}^\prime$ in $M^\prime$. Given a
leaf $L\subset M$ of ${\mathcal F}$, any point $p\in L$ and a small
(holomorphic) transverse section $\Sigma\subset M$, with $p \in
\Sigma \cap L$, $\Sigma$ transverse to $L$ and ${\mathcal F}$,  we
may introduce the {\it holonomy group} $\Hol({\mathcal
F},L,\Sigma,p)$ of this leaf of ${\mathcal F}$ as the usual holonomy
group of the leaf $L\in {\mathcal F}^\prime$ of the foliation
${\mathcal F}^\prime$ calculated with respect to the transverse
section $\Sigma\subset M^\prime$ at the point $p\in M^\prime$ (see
\cite{Camacho-LinsNeto} for instance).

 In case $M$ is a projective manifold (for instance, the $n+1$ dimensional  complex projective space
$M=\mathbb CP^{n+1}$), the foliation ${\mathcal F}$ is necessarily
{\it algebraic} in the sense that it is given by algebraic equations
(polynomial vector fields) in any affine subspace of $M$. In this
case by an {\it algebraic solution} of ${\mathcal F}$ we mean an
algebraic curve $\Lambda \subset M$ such that $M\setminus (\Lambda
\cap \sing({\mathcal F}))$ is a union of leaves of ${\mathcal F}$.

\begin{Proposition}
\label{Proposition:holonomygroups}
Let ${\mathcal F}$ be a dimension one holomorphic  foliation with
singularities of the complex projective space $\mathbb C P^{n+1}$.
Let $\Lambda\subset \mathbb C P^{n+1}$ be an algebraic solution of
${\mathcal F}$. Then the holonomy groups of the leaves $L\subset
\Lambda$ are contained in algebraic groups.

\end{Proposition}

\begin{proof}
Indeed, given a point $p \in \Lambda\setminus \sing({\mathcal F})$
we may choose a liner hyperplane $\mathbb E(n)\subset \mathbb
CP^{n+1}$ such that $\Lambda$ and $\mathbb E(n)$ meet transversely at
$p$. The choice of an affine system of coordinates $\mathbb C^{n+1}
_p $ centered at $p$ gives a polynomial vector field $X_p$ with
isolated singularities that defines the foliation ${\mathcal
F}\big|_{\mathbb C^{n+1}_p}$ in the ordinary sense: the leaves of
${\mathcal F}\big|_{\mathbb C^{n+1}_p}$ are the non-singular
integral curves of $X_p$ in $\mathbb C^{n+1}_p$. Then we choose $n$
polynomial one-forms $\omega_1,...,\omega_n$ in the affine space
$\mathbb C^{n+1}_p$ with the property that $\omega_j \cdot X=1$ and
such that the exterior product $\omega_1 \wedge \cdots \wedge
\omega_n \not\equiv 0$. Denote by $\omega_j ^*$ the restriction of
$\omega_j$ to $\mathbb E(n) \cap \mathbb C^{n+1}_p \simeq \mathbb
C^n$. Then each $\omega_j ^*$ is rational. If we denote by $\Sigma$
the germ at $p$ of disc induced by $\mathbb E(n)$, then  for each
holonomy diffeomorphism $g \in \Hol({\mathcal F}, \Lambda,
\Sigma,p)$ we have $g^*(\omega_j)\wedge \omega_1\wedge \cdots
\omega_n\equiv 0$ because the holonomy maps preserve the leaves of
the foliation. Therefore, the holonomy group of the leaves contained
in $\Lambda$ are contained in algebraic groups.
\end{proof}

\section{Abelian groups tangent to the identity}
\label{section:flatabelianclassification}

Now we study the characterization and classification of  abelian
tangent to the identity groups,  which is the subject of
Theorem~\ref{Theorem:A}.

According to \cite{Brochero}  (Proposition 4.1) {\it  every nilpotent subalgebra
$\mathcal L$  of $\hat{\mathcal X}(\mathbb C^2, 0)$ is metabelian}.
Also following \cite{Brochero}, this proposition implies the following characterization of abelian
subgroups of  $\Diff_{\Id}(\mathbb C^2, 0)$.

\begin{Proposition}[cf. \cite{Brochero}
Corollary 4.4]
\label{Proposition:flatabeliananalytic}
 If $ G < \Diff_{\Id}(\mathbb C^2, 0)$ is a  abelian tangent to the identity
 {\rm(}convergent{\rm)} group, then
one of the following items is true: \begin{itemize}

\item[{\rm(1)}] There is a formal vector field $\hat {X}$
with $\exp (\hat {X}) \in  G$ such that for each $g\in G$ we have
$g=\exp(T_{g} \hat {X})$  where $T_{g}\in K_2$ is a rational function such that $\hat {X}(T_{g}) = 0$;

\item[{\rm(2)}] $ G < \left<{f}^{[t]} \circ { g}^{[s]}| t,s \in \mathbb  C\right>$,
where $f,  g\in  G$ and $[ f,  g] =\Id$.

\end{itemize}
\end{Proposition}

Notice that Proposition~\ref{Proposition:flatabeliananalytic} above
is for convergent (analytic) objects. In this paper we obtain an extension of this result, Theorem~\ref{Theorem:nilpotent} below. This is  based also in the following  formal version of Proposition~\ref{Proposition:flatabeliananalytic} which is easily
obtained by a  mimic of its proof:

\begin{Proposition}
\label{proposition:abeliantangentidentity}
 Let $ G < {\hat
\Diff_{\Id}}(\mathbb{C}^2,0)$ be an abelian  subgroup tangent to the
identity, we have the two following possibilities:

\begin{enumerate}
\item There is a formal vector field $\hat {X}$, invariant by $ G$, such that
for each element $\hat f\in  G$ there is a rational function $T\in K_2$
depending on $\hat f$, such that $\hat {X}(N)=0$ and $\hat
f=\exp(N\hat {X})$.

\item There are formal commuting vector fields $\hat {X}$ and $\hat {Y}$
such that   $\exp(\hat {X}),\exp(\hat {Y})\in  G$ and $G < \langle
\exp(t\hat {X})\circ \exp(s\hat {Y})\mid t,s\in\mathbb{C}\rangle$.

\end{enumerate}
\end{Proposition}

\begin{Definition}[Formal meromorphic one-forms]
\label{Definition:formalobjects} {\rm By a {\it formal meromorphic
function} of $n$ complex variables we shall mean a formal quotient
$\hat R=\frac{\hat P}{\hat Q}$  of two formal power series with
positive exponents $\hat P, \hat Q \in \hat{\mathcal O_n}=\mathbb
C[[z]]$. In other words, the field of formal meromorphic functions
of $n$ variables $\hat{K}_{n}$ will be the fraction field of the
domain of integrity $\hat{\mathcal O_n}$. By a {\it meromorphic}
formal one-form we  mean a formal expression $\hat \omega=
\sum\limits_{j=1}^n \hat R_j \, dz_j$ where each $\hat R_j$ is a
formal meromorphic function as defined above. The exterior
derivative, wedge product and other concepts are defined for
meromorphic formal one-forms in the same way as for analytic
one-forms.}
\end{Definition}

\begin{proof}[Proof of Theorem~\ref{Theorem:A}]
Assume that the Lie algebra $l( G_{\Id})$ has dimension one. By
Proposition~\ref{proposition:abeliantangentidentity}  there is a
formal vector field $\hat {X}$, invariant by $ G$, such that for
each $\hat f\in  G$ there is a rational function $T=T_{\hat f}\in K_2$ with
$\hat {X}(T)=0$ and $\hat f=\exp(T\hat {X})$. Suppose that for some
$\hat f\in G$ the function  $T=T_{\hat f}$ is not constant. We
consider the one-form $\hat \omega:=dT$. This is a non-trivial
closed rational one-form and we claim that this is $ G$-invariant.
In fact, take $\hat g\in G$ and write $G=\exp(S \hat {X})$ for some
rational function  $S\in K_2$ such that $\hat {X}(S)=0$. Then $(S\hat
{X})(T)=dT(S\hat {X})=S dT(\hat {X})=0$.  Therefore $T\circ
\exp(S\hat {X})=T$. This gives $\hat g^*(\omega)=\hat
g^*(dT)=d(T\circ \hat g)=d(T\circ \exp(S\hat {X}))=dT=\omega$,
proving the claim.  This corresponds to (i) in
Theorem~\ref{Theorem:A}. \vglue.2in

Now we consider the where $T_{\hat f}$ is constant for each $\hat
f\in  G$. In this case each element $\hat f\in  G$ writes as $\hat
f=\exp(c_{\hat f} \hat {X})$ for some constant $c_{\hat f}\in
\mathbb C$. In other words, $ G$ embeds into the flow of $\hat {X}$
as in (ii) in the statement.

Suppose now that $ G$ is as in (2) in
Proposition~\ref{proposition:abeliantangentidentity}. Then there are
two $\mathbb C$-linearly independent, formal, commuting vector
fields $\hat {X}_j$, invariant by $G$,  such that $\exp(\hat
{X}_j)\in  G$ and $G < \langle \exp(t\hat {X}_1)\circ \exp(s\hat
{X}_2)\mid t,s\in\mathbb{C}\rangle$.
\end{proof}

As already mentioned in the Introduction, (iii) in Theorem~\ref{Theorem:A} admits a converse,
proved as follows:

\begin{proof}[Proof of Theorem~\ref{Theorem:B} - First Part]
By hypothesis $G< \hat{\Diff}(\mathbb C^n,0)$  admits two linearly independent invariant formal vector fields $\hat X_1, \hat X_2$. We fix coordinates $(x,y)$ and  write $\hat {X}_j =
A_j\frac{\partial}{\partial x} + B_j\frac{\partial}{\partial y}$.
Since   $\hat {X}_j$ is $G$-invariant, we have
\[
\begin{bmatrix}
\frac{\partial g_1}{\partial x} & \frac{\partial g_1}{\partial y}  \\[7pt]
\frac{\partial g_2}{\partial x} & \frac{\partial g_2}{\partial y}
\end{bmatrix}
\begin{bmatrix}
A_1(z) & A_2(z) \\[7pt] B_1(z) & B_2(z)
\end{bmatrix}
=
\begin{bmatrix}
A_1(\hat g) & A_2(\hat g) \\[7pt] B_1(\hat g) & B_2(\hat g)
\end{bmatrix}
\]
Taking transposes, we obtain:
\[
\begin{bmatrix}
A_1(z) & B_1(z) \\[7pt] A_2(z) & B_2(z)
\end{bmatrix}
\begin{bmatrix}
\frac{\partial g_1}{\partial x} & \frac{\partial g_2}{\partial x}  \\[7pt]
\frac{\partial g_1}{\partial y} & \frac{\partial g_2}{\partial y}
\end{bmatrix}
=
\begin{bmatrix}
A_1(\hat g) & B_1(\hat g) \\[7pt] A_2(\hat g) & B_2(\hat g)
\end{bmatrix}
\]
thus
\begin{equation}
\label{equation:invar}
\begin{bmatrix}
\frac{\partial g_1}{\partial x} & \frac{\partial g_2}{\partial x}  \\[7pt]
\frac{\partial g_1}{\partial y} & \frac{\partial g_2}{\partial y}
\end{bmatrix}
\begin{bmatrix}
A_1(\hat g) & B_1(\hat g) \\[7pt] A_2(\hat g) & B_2(\hat g)
\end{bmatrix}^{-1}
=
\begin{bmatrix}
A_1(z) & B_1(z) \\[7pt] A_2(z) & B_2(z)
\end{bmatrix}^{-1}
\end{equation}
Based on this last equation we  take:
\begin{equation}
\label{equation:closedform}
\begin{bmatrix}
C_1(z) & C_2(z) \\[7pt] D_1(z) & D_2(z)
\end{bmatrix}
=
\begin{bmatrix}
A_1(z) & B_1(z) \\[7pt] A_2(z) & B_2(z)
\end{bmatrix}^{-1}
= \frac{1}{Q(z)}
\begin{bmatrix}
B_2(z) & -B_1(z) \\[7pt] -A_2(z) & A_1(z)
\end{bmatrix}
\end{equation}

Where $Q(z) = A_1B_2 - A_2B_1$. Now we define $\hat \omega_j := C_jdx +
D_jdy$. By the above equations~\eqref{equation:invar} and \eqref{equation:closedform}, the one-forms  $\hat \omega_j$ are invariant
for $ G$, i.e, for each
 $\hat
g\in G$, we have  $\hat g^*(\hat \omega_j) = \hat \omega_j$ $(j = 1,2)$.
These forms are  $\mathbb C$-linearly independent (cf. \eqref{equation:closedform}). Let us now  show that the $\hat \omega_j$ are closed forms. We have to show that
$\frac{\partial D_j}{\partial x} - \frac{\partial C_j}{\partial y} =
0$. Since $[\hat X_1,\hat {X}_2] = 0$ then
\begin{eqnarray*}
\frac{\partial A_2}{\partial x}A_1 + \frac{\partial A_2}{\partial y}B_1 &
=  & \frac{\partial A_1}{\partial x}A_2 + \frac{\partial A_1}{\partial y}B_2\\
\frac{\partial B_2}{\partial x}A_1 + \frac{\partial B_2}{\partial y}B_1 &
= & \frac{\partial B_1}{\partial x}A_2 + \frac{\partial B_1}{\partial y}B_2
\end{eqnarray*}
Thus
\begin{align*}
& Q^2(\frac{\partial D_1}{\partial x} - \frac{\partial C_1}{\partial
y})= Q^2(-\frac{1}{Q}\frac{\partial A_2}{\partial x} +
\frac{A_2}{Q^2}\frac{\partial Q}{\partial x} -
\frac{1}{Q}\frac{\partial B_2}{\partial y} +
\frac{B_2}{Q^2}\frac{\partial Q}{\partial y})\\
&= A_2\frac{\partial Q}{\partial x} - Q\frac{\partial A_2}{\partial x} +
B_2\frac{\partial Q}{\partial y} - Q\frac{\partial B_2}{\partial y}\\
&= A_2B_2\frac{\partial A_1}{\partial x} + A_2A_1\frac{\partial B_2}{\partial x} -
A_2B_1\frac{\partial A_2}{\partial x} - A_2^2\frac{\partial B_1}{\partial x} -
A_1B_2\frac{\partial A_2}{\partial x} + A_2B_1\frac{\partial A_2}{\partial x} +\\
& \hspace{0.5cm}
B_2^2\frac{\partial A_1}{\partial y} + B_2A_1\frac{\partial B_2}{\partial y} -
B_2B_1\frac{\partial A_2}{\partial y} - B_2A_2\frac{\partial B_1}{\partial y} -
A_1B_2\frac{\partial B_2}{\partial y} + A_2B_1\frac{\partial B_2}{\partial y}\\
&= B_2(\frac{\partial A_1}{\partial x}A_2 + \frac{\partial A_1}{\partial y}B_2 -
\frac{\partial A_2}{\partial x}A_1 - \frac{\partial A_2}{\partial y}B_1) +\\
& \hspace{0.5cm}
A_2(\frac{\partial B_2}{\partial x}A_1 + \frac{\partial B_2}{\partial y}B_1 -
\frac{\partial B_1}{\partial x}A_2 - \frac{\partial B_1}{\partial y}B_2)\\
&= B_2.0 + A_2.0 = 0.
\end{align*}

Therefore $\hat \omega_1$ is closed. Similarly $\hat \omega_2$ is a closed one-form.  This proves the first part of Theorem~\ref{Theorem:B}.  Assume now that $G$ is tangent to the identity.
Write $\hat \omega_j = A_j dx + B_j dy$ then,
\[
\hat \omega_j = \hat g^*(\hat \omega_j) = A_j(\hat g) dg_1 + B_j(\hat g) dg_2
\]
\[
=(A_j(\hat g)\frac{\partial g_1}{\partial x} + B_j(\hat g)\frac{\partial
g_2}{\partial x}) dx + (A_j(\hat g)\frac{\partial g_1}{\partial y} +
B_j(\hat g)\frac{\partial g_2}{\partial y}) dy
\]
Thus
\[
\begin{bmatrix}
A_1(\hat g) & B_1(\hat g) \\[7pt] A_2(\hat g) & B_2(\hat g)
\end{bmatrix}
\begin{bmatrix}
\frac{\partial g_1}{\partial x} & \frac{\partial g_2}{\partial x}  \\[7pt]
\frac{\partial g_1}{\partial y} & \frac{\partial g_2}{\partial y}
\end{bmatrix}
=
\begin{bmatrix}
A_1(z) & B_1(z) \\[7pt] A_2(z) & B_2(z)
\end{bmatrix}
\]
consequently
\[
\begin{bmatrix}
\frac{\partial g_1}{\partial x} & \frac{\partial g_2}{\partial x}  \\[7pt]
\frac{\partial g_1}{\partial y} & \frac{\partial g_2}{\partial y}
\end{bmatrix}
\begin{bmatrix}
A_1(z) & B_1(z) \\[7pt] A_2(z) & B_2(z)
\end{bmatrix}^{-1}
=
\begin{bmatrix}
A_1(\hat g) & B_1(\hat g) \\[7pt] A_2(\hat g) & B_2(\hat g)
\end{bmatrix}^{-1}
\]
Let us introduce $\hat {X}_1, \hat {X}_2$ as follows:
\[
\hat {X}_1 = \frac{1}{Q(z)}(B_2\frac{\partial}{\partial x} -
A_2\frac{\partial}{\partial y}) \, ,\hspace{0.3cm}\, \, \, \,
\hspace{0.3cm} \hat {X}_2 =
\frac{1}{Q(z)}(-B_1\frac{\partial}{\partial x} +
A_1\frac{\partial}{\partial y})
\]
where $Q(z) = A_1(z)B_2(z) - A_2(z)B_1(z)$. Since the $\hat
\omega_j$ are closed one-forms we have $\hat g^*(\hat {X}_j) = \hat
{X}_j$,  for all $\hat g\in G$ and also $[\hat {X}_1, \hat {X}_2] = 0$.
Also  $\hat {X}_1, \hat {X}_2$ are $\mathbb C$-linearly independent
in $\hat{K}_2$.
Note that $\{\hat {X}_1, \hat {X}_2\}$ is a basis for the vector
space $\hat{\mathcal X}(\mathbb{C}^{2},0)\otimes \hat{K}_{2}$ and,
since for $\hat g\in G$ we can write $\hat g = \exp(\hat {Y}_{\hat
g})$  with $\hat {Y}_{\hat g}\in \hat{\mathcal
X}(\mathbb{C}^{2},0)$,
then $ \hat {Y}_{\hat g} = u_1\hat {X}_1 + u_2\hat {X}_2$, where
$u_j\in\hat{K}_2$. On the other side $\hat g^*(\hat {X}_j) = \hat
{X}_j$ we have that $[\hat {Y}_{\hat g},\hat {X}_1] = 0$ and $[\hat
{Y}_{\hat g},\hat {X}_2] = 0$, then $\hat {X}_j(u_k) = 0$ $(j,k =
1,2)$, consequently $u_j$ are constant in $C^*$. Now if $G < \langle
\exp(t\hat {Y}_{\hat g}))\mid t\in\mathbb{C}\rangle$ there is
nothing left to prove, thus suppose that there is $\hat h\in G$,
$\hat h = \exp(\hat {Y}_{\hat h})$ with $\hat {Y}_{\hat g}$ and
$\hat {Y}_{\hat h}$ $\mathbb C$-linearly independent in $\hat{K}_2$
then there are $v_j\in\mathbb{C}$ such that $ \hat {Y}_{\hat h} =
v_1\hat {X}_1 + v_2\hat {X}_2$ and $(u_1,u_2)$,$(v_1,v_2)$ are
$\mathbb C$-linearly independent in $\mathbb{C}^{2}$, therefore
$\hat {X}_j\in\hat{\mathcal X}(\mathbb{C}^{2},0)$ and $G = \exp(
a\hat {X}_1)\circ\exp (b\hat {X}_2)$, with $a,b\in C^*$, therefore
there are formal vector fields $\hat {X}$ and $\hat {Y}$ such that
$[\hat {X},\hat {Y}] = 0$ and $G < \langle \exp(t\hat
{X})\circ\exp(s\hat {Y})\mid t,s\in\mathbb{C}\rangle$. Therefore $
G$ is abelian.  This proves the first part of
Theorem~\ref{Theorem:B}.
\end{proof}

The second part will be concluded in the next section (cf. Proposition~\ref{proposition:converseabelian}).

\section{Groups preserving closed one-forms}
\label{section:integrationlemma} In this section we
finish the proof of Theorem~\ref{Theorem:B}. Indeed, we proceed
studying the classification of groups of formal diffeomorphisms
preserving closed meromorphic one-forms in $(\mathbb{C}^2,0)$.
Special attention is given to the ``generic" case where the group
exhibits two transverse formal separatrices. Before going further
into the main subject we recall some classical facts about
integration of closed meromorphic one-forms in several complex
variables.

\begin{Proposition}[Integration Lemma, cf. \cite{CLNS}, \cite{Scardua} Example 1.6]
\label{Proposition:integrationlemma} Let $\omega$ be a closed
meromorphic  one-form on $M$ where $M$ is
 a polydisc in $\mathbb C^n$. Then there are irreducible holomorphic functions
 $f_1,..,f_r\in \mathcal O(M)$, $n_1,...,n_r \in \mathbb N$, complex
  numbers $\lambda_1,...,\lambda_r$ and a
 holomorphic function $g \in \mathcal O(M)$ such that
 \[
 \omega=\sum\limits_{j=1}^r \lambda_j \frac{df_j}{f_j} +
 d(\frac{g}{f_1^{n_1}\cdots f_r ^{n_r}})
 \]
 \end{Proposition}

The polar set of $\omega$ is given in irreducible components by
$\bigcup\limits_{j=1}^r \{f_j=0\}$, $n_j$ is the order of
$\{f_j=0\}$ as a component of the polar set of $\omega$, $\lambda_j$
is the residue of $\omega$ at the component $\{f_j=0\}$ and the
function $g$ has no common factors with $f_j$ in $\mathcal O(M)$.

If $M=\mathbb C^n$   and $\omega$ is a rational closed one-form then we have the same
result, where the $f_j$ are irreducible polynomials and $g$ is a
polynomial without common factors with the $f_j$. The proof of
Theorem~\ref{Proposition:integrationlemma} relies on integration and
the fact that the first homology group of the complement of a pure
codimension one analytic subset $\Lambda=\bigcup\limits_{j=1}^r
\Lambda_j$, where each $\Lambda_j$ is an irreducible component, of a
polydisc $M$ as above, is generated by small loops around the
components $\Lambda_j$, contained in transverse discs circulating
the component. Then a standard argument involving Laurent series
implies the result. This cannot be repeated in the formal case,
because we cannot rely on integration processes, at  first glance.
Nevertheless, we still have a formal  version of
Theorem~\ref{Proposition:integrationlemma} as follows:

\begin{Proposition}[Formal integration lemma]
\label{Proposition:formalintegrationlemma} Let $\hat\omega$ be  a
closed formal meromorphic one-form in $n$ complex variables. Denote
by $\hat f_j\in\hat {\mathcal O}_k, \, j=1,...,r$ the formal
equations of the set of poles of $\hat \omega$, in independent
terms. Then, there are $\lambda_j \in \mathbb C$ and $n_j \in
\mathbb N$ and a formal function $\hat g \in \hat{\mathcal O}_k$
such that
\[
\hat \omega= \sum\limits_{j=1} ^r \lambda_j \frac{d\hat f_j}{\hat
f_j} + d(\frac{\hat g}{\hat f_1 ^{n_1}\cdots \hat f_j^{n_r}})
\]
\end{Proposition}
The proof is somehow  similar to the proof of the local analytic
version and it is based on the following:

\begin{Lemma}
A closed formal meromorphic one-form $\hat \omega$ in $n$ complex
variables, without residues
 is exact: $\hat \omega=d \hat f$ for some meromorphic
formal function $\hat f \in \hat{K}_{n}$.
\end{Lemma}
This lemma is proved analogously to the following particular case:

\begin{Lemma}
\label{Lemma:integrationdimensiontworesiduezero} Let $\hat \omega$
be a closed formal meromorphic one-form in   two complex variables
and assume that the polar set of $\hat \omega$ consists of two
transverse formal curves, and that the residues of $\hat \omega$ are
all zero. Then $\hat \omega$ is exact, indeed, in suitable formal
coordinates $(x,y)$ we can write
\[
\hat \omega= d(\frac{\hat f}{x^n y ^m})
\]
for some $n, m \in \mathbb N$ and some formal function  $\hat f \in
\hat{\mathcal O}_2$.
\end{Lemma}

\begin{proof}
Since the polar set of $\hat \omega$ consists of two transverse
formal curves, we can find formal coordinates $(x,y)$ such that this
polar set corresponds to the coordinates axes. We write $\hat \omega
= (\hat P/x^{n+1} y^{m+1}) dx +  (\hat Q/x^{n+1} y^{m+1})dy$ where
$\hat P ,  \, \hat Q \in \mathbb C[[x,y]]$. We can write $\hat P
=\sum\limits_{\nu=0} ^\infty P_\nu$ and $\hat Q = \sum
\limits_{\nu=0} ^\infty Q_\nu$ in terms of homogeneous polynomials
$P_\nu, \, Q_\nu$ of degree $\nu-n-m$. Then $\hat \omega=
\sum\limits_{\nu = 0 } ^\infty ( P_\nu/x^{n+1} y^{m+1}) dx +
(Q_\nu/x^{n+1} y^{m+1})dy =\sum\limits_{\nu=-n -m} ^\infty
\omega_{\nu}$ where $\omega_\nu=  (P_\nu/x^{n+1} y^{m+1}) dx +
(Q_\nu/x^{n+1} y^{m+1})dy$ is a homogeneous rational one-form of
degree $\nu - n -m -2$. Then $d \hat \omega = \sum\limits_{\nu=-n
-m-2} ^\infty d\omega_{\nu}$ where each one-form $d \omega_\nu$ is
homogeneous of degree $\nu-1$. Therefore, since $\hat \omega$ is
closed  we have $0 = d \hat \omega= \sum\limits_{\nu=-n-m}^\infty
d\omega_\nu$ and then $d\omega_\nu=0, \, \forall \nu \geq -n -m$.
Since  $\hat \omega$ has no residues, the same holds for
$\omega_\nu$. Moreover, because each form $\omega_\nu$ is of the
form  $ \omega_\nu=P_\nu/x^{n+1} y^{m+1}) dx +  (Q_\nu/x^{n+1}
y^{m+1})dy$, we conclude from the Integration lemma that $\omega_\nu
= d(\frac{ f_\nu}{x^{n}y^{m}})$ for some homogeneous polynomial
$f_\nu$ of degree $\nu$. Thus $\hat \omega=d(
{\sum\limits_{\nu=0}^\infty f_\nu}/{x^{n}y^{m}})=d(\hat f/ x^n y^m)$
where $\hat f= \sum\limits_{\nu=0}^\infty f_\nu\in \hat{\mathcal
O}_2$.

\end{proof}

As a consequence we obtain the following particular  case of
Proposition~\ref{Proposition:formalintegrationlemma}:

\begin{Proposition}
\label{Proposition:integrationlemmadimtwo} Let $\hat \omega$ be a
closed formal meromorphic one-form in two complex variables and
assume that the polar set of $\hat \omega$ consists of  two
transverse formal curves. Then $\hat \omega$ writes  in suitable
formal coordinates $(x,y)$ as
\[
\hat \omega= \lambda \frac{dx}{x} + \mu \frac{dy}{y} + d(\frac{\hat f}{x^n y ^m})
\]
for some $\lambda, \, \mu \in \mathbb C$, some $n, m \in \mathbb N$
and some formal function $\hat f \in \hat{\mathcal O}_2$.
\end{Proposition}

\begin{proof} As in the proof of
Lemma~\ref{Lemma:integrationdimensiontworesiduezero} we choose
formal coordinates $(x,y)$ such that the polar set of $\hat \omega$
corresponds to the coordinate axes. Denote by $\lambda\in\mathbb C$
and $\mu\in \mathbb C$ the residue of $\hat \omega$ at the $x$-axis
and $y$-axis respectively. Then $\hat \theta=\hat \omega - \lambda
\frac{dx}{x} + \mu \frac{dy}{y}$ is a closed formal meromorphic
one-form with polar set contained in the coordinate axes and zero
residues. By Lemma~\ref{Lemma:integrationdimensiontworesiduezero}
we can write $\hat \theta= d(\frac{\hat f}{x^n y ^m})$ for some
formal function $\hat f \in \hat{\mathcal O}_2$.
\end{proof}

An improvement of the above proposition is the following:

\begin{Lemma}
\label{clasification} Let $\hat \omega_j, \, j=1,2$ be  $\mathbb
C$-linearly independent closed formal meromorphic one-forms in two
variables with polar sets along two transverse formal curves. Then
there are formal coordinates $(x,y)$  such that each $\hat \omega_j$
writes:
\begin{equation}
\label{equation:integrationlemma}
 \hat \omega_j = a_j\frac{dx}{x} +
b_j\frac{dy}{y} + d(\frac{c_j}{x^{n_j}y^{m_j}})
\end{equation}
for some constant $a_j, b_j, c_j \in \mathbb C$ and some $n_j, m_j
\in \mathbb N$.
\end{Lemma}

\begin{proof}
By Proposition~\ref{Proposition:integrationlemmadimtwo} we can write
 $\hat\omega_j =
a_j\frac{dx}{x} + b_j\frac{dy}{y} + d(\frac{\hat
f_j}{x^{n_j}y^{m_j}})$, where $a_j, b_j\in\mathbb{C}; \, n_j, m_j
\in \mathbb N$ and $\hat f_j\in \hat{\mathcal O}_k$. Let us write
$n_1=n, m_1=m$ and $n_2=p, m_2 =q$.  We take a model of formal
change of coordinates $\hat \phi = (xu,yv)$, where we want that
$\hat\phi^*(a_1\frac{dx}{x} + b_1\frac{dy}{y} +
d(\frac{c_1}{x^ny^m})) = \hat \omega_1$ and
$\hat\phi^*(a_2\frac{dx}{x} + b_2\frac{dy}{y} +
d(\frac{c_2}{x^py^q})) = \hat \omega_2$ implies $a_1\frac{du}{u} +
b_1\frac{dv}{v} = d(\frac{1}{x^ny^m}.(\hat f_1 -
\frac{c_1}{u^nv^m}))$ and $a_2\frac{du}{u} + b_2\frac{dv}{v} =
d(\frac{1}{x^py^q}.(\hat f_2 - \frac{c_2}{u^pv^q}))$ then,
\[
(a_1\ln u +b_1\ln v)x^ny^m - \hat f_1 + \frac{c_1}{u^nv^m} +
k_1x^ny^m = 0
\]

and \[
 (a_2\ln u +b_2\ln
v)x^py^q - \hat f_2 + \frac{c_2}{u^pv^q} + k_2x^py^q = 0
\]
Now define a formal meromorphic function $\hat R$ by
\[
\hat R(x,y,u,v) = (\hat R_1(x,y,u,v), \hat R_2 (x,y,u,v)) \] where
\[
\hat R_1(x,y,u,v)= (a_1\ln u +b_1\ln v)x^ny^m -
 \hat f_1 + \frac{c_1}{u^nv^m} + k_1x^ny^m
 \]
 and
 \[
 \hat R_2(x,y,u,v)= (a_2\ln u +b_2\ln v)x^py^q - \hat f_2 + \frac{c_2}{u^pv^q} +
 k_2x^py^q.
 \]
We have $\hat R(0,0,u,v) = (0,0)$, so that if $c_1 = \hat f_1(0)$
and $c_2 = \hat f_2(0)$ we have $\frac{1}{u^n(0)v^m(0)} = 1$  and
$\frac{1}{u^p(0)v^q(0)} = 1$ and as $\Det (J_2R(0,(u,v)) = (nq -
mp)u^{k + n + 1}v^{m + q + 1} \neq 0$, if $(m,n)$ and $(p,q)$ are
$\mathbb C$-linearly independent,  from the formal version of the
Implicit function theorem we obtain a unique solution $(u,v)$.
\end{proof}

\begin{Remark}
{\rm Let  $ G < {\hat \Diff_{\Id}}(\mathbb{C}^2,0)$ be a
subgroup. Given a closed meromorphic 1-form  $\hat \omega$ such that
$\hat \omega$ is invariant by $ G$, if $\hat \omega$ is conjugated
to a 1-form $\hat\alpha$ by a diffeomorphism $\hat h$, then the
1-form $\hat\alpha$ is invariant by the group $\hat h^{-1}\circ
\hat g \circ \hat h$. As the groups $ G$ and $\hat h^{-1}\circ  \hat
g \circ \hat h$ have similar algebraic proprieties, there is no loss
of generality in assuming that the forms are as in the normal form
of Lemma~\ref{clasification}. }
\end{Remark}

The second part of the proof of Theorem~\ref{Theorem:B}  follows
from the following proposition:
\begin{Proposition}
\label{proposition:converseabelian}
 Let $ G < \hat {\Diff}(\mathbb{C}^2,0)$ be a
 subgroup with two transverse separatrices, if there are
 $\mathbb C$-linearly independent  closed formal meromorphic formal  one-forms $\hat \omega_j, \,
(j=1,2)$ which are invariant by $ G$ then $ G$ is abelian. Indeed,
 $ G$ is formally conjugate to a group of  diffeomorphisms
 generated by any of the following types:

\begin{enumerate}
\item[{\rm(a)}]
$ \hat g(x,y) = (x\frac{(1 +
\frac{k_2}{c_2}x^py^q)^{\frac{m}{D}}}{(1 +
\frac{k_1}{c_1}x^ny^m)^{\frac{q}{D}}},
 y\frac{(1 + \frac{k_1}{c_1}x^ny^m)^{\frac{p}{D}}}{(1
  + \frac{k_2}{c_2}x^py^q)^{\frac{n}{D}}})
$
\item[{\rm(b)}]
$\hat g(x,y) = (ax, \frac{a^{-\frac{n}{m}}y}{(1 +
kx^ny^m)^{\frac{1}{m}}})$
\item[{\rm(c)}]
$\hat g(x,y) = (\frac{b^{-\frac{m}{n}}x}{(1 +
kx^ny^m)^{\frac{1}{n}}}, by)$
\end{enumerate}

\end{Proposition}
\begin{proof}
 A diffeomorphism $\hat g \in
G$ writes $\hat g(x,y) = (xu,yv)$ where $u, v \in \hat {\mathcal
O}_2$ satisfy $\hat u(0)\ne 0, \hat v(0)\ne 0$. From
equation~\eqref{equation:integrationlemma} and $\hat g^*(\hat
\omega_j) = \hat \omega_j$ we obtain:
\[
a_1\ln u + b_1\ln v = \frac{c_1}{x^ny^m}.(1 - \frac{1}{u^nv^m}) +
k_1 {\rm \hspace{0.3cm} and} \hspace{0.3cm}
 a_2\ln u + b_2\ln v = \frac{c_2}{x^py^q}.(1 - \frac{1}{u^pv^q}) + k_2
\]

If $a_1 = a_2 = b_1 = b_2 = 0$ then $u^nv^m = \frac{1}{1 +
\frac{k_1}{c_1}x^ny^m}$ and $u^pv^q =
 \frac{1}{1 + \frac{k_2}{c_2}x^py^q}$, as $(m,n)$ and $(p,q)$
 must be $\mathbb C$-linearly independent we have
\[
\hat g(x,y) =
(x\frac{(1 + \frac{k_2}{c_2}x^py^q)^{\frac{m}{D}}}{(1
 + \frac{k_1}{c_1}x^ny^m)^{\frac{q}{D}}},
 y\frac{(1 + \frac{k_1}{c_1}x^ny^m)^{\frac{p}{D}}}{(1
  + \frac{k_2}{c_2}x^py^q)^{\frac{n}{D}}})
\]
with $D = nq - pm$. The group $ G$ therefore has just linear
diffeomorphisms as above and  is an abelian group.

Assume now that  the left side of equality is holomorphic we have,
$\frac{1}{u^nv^m} = 1$ and $\frac{1}{u^pv^q} = 1$ when $(m,n)$ and
$(p,q)$ are $\mathbb C$-linearly independent, we have that $u$ and
$v$ are constant, so that $G$ is linear therefore $ G$ is abelian.

A similar argumentation with the other possible cases gives the
forms:
\[
\hat g(x,y) = (ax, \frac{a^{-\frac{n}{m}}y}{(1 +
kx^ny^m)^{\frac{1}{m}}}), \hat g(x,y) =
(\frac{b^{-\frac{m}{n}}x}{(1 + kx^ny^m)^{\frac{1}{n}}}, by)
\]
 and
\[
\hat g(x,y) = (\frac{x}{(1 + k_1x^n)^{\frac{1}{n}}}, \frac{y}{(1 +
k_2y^m)^{\frac{1}{m}}}).
\]
 In particular, on each case,  $ G $ is
abelian.
\end{proof}

\begin{Remark}[Holomorphic case]
\label{Remark:withoutpoles}{\rm  If $ G$ is invariant by two
$\mathbb C$-linearly independent closed formal one-forms (without
poles) then $ G = \{\Id\}$.}
\end{Remark}
\begin{proof}
Let $\hat g\in G$ and write  $\hat \omega_j=d\hat f_j$, $( j = 1,
2)$ left invariant by $ G$. Take $\hat \Phi = (\hat f_1,\hat f_2)$,
as $d\hat f_1$ and $d\hat f_2$ are $\mathbb C$-linearly independent
in this neighborhood of the origin, we have that $\hat \Phi$ is a
formal diffeomorphism. Therefore we may assume that
$\hat\omega_1=dx$ and $\hat\omega_2=dy$, {\it i.e.},  $(\hat
\Phi\circ \hat g \circ \hat \Phi^{-1})^*(dx) = dx$ and $(\hat
\Phi\circ \hat g \circ \hat \Phi^{-1})^*(dy) = dy$, because $(\hat
\Phi\circ \hat g \circ \hat \Phi^{-1})^*(dx) = (\hat
\Phi^{-1})^*\circ \hat g^*\circ\hat \Phi^*(dx) = (\hat
\Phi^{-1})^*\circ \hat g^*(d\hat f_1) = (\hat \Phi^{-1})^*(d\hat
f_1) = (\hat \Phi^{-1})^*\circ\hat \Phi^*(dx) = dx$, therefore $\hat
\Phi\circ \hat g \circ \hat \Phi^{-1} = \{\Id\}$
 and consequently $ G = \{\Id\}$.
\end{proof}

\section{Metabelian groups}
\label{section:metabelian} Now we  study metabelian groups in
$\hat\Diff(\mathbb C^2, 0)$, that is, subgroups $G< \hat
\Diff(\mathbb C^2, 0)$  such that the subgroup of commutators
$G^{(1)}=[G,G]$ is abelian. Let $G$ be such a metabelian subgroup.
Then, the derivative group  $DG< \GL(\mathbb C, 2)$ is also
metabelian but not necessarily  abelian. For instance, take $G$ as
the linear subgroup of $2 \times2$ triangular superior matrices.
Then  $G$ is not abelian but $G^{(1)}$  is abelian.

Now if the group  $D G$ is abelian then $G^{(1)}$ is  tangent to the
identity, which is a very useful property. For this reason, in our
statements related to this case (Theorem~\ref{Theorem:metabelian}), we may alternatively require that $D G$ is abelian.

\begin{Lemma}
Let $G < \hat{\Diff}(\mathbb{C}^{2},0)$ be a subgroup with $D G$
abelian. Suppose that, there are two $\mathbb C$-linearly independent
vector fields $\hat {X}$ and $\hat {Y}$,
projectively invariant by $ G$ and such that $[\hat {X},\hat {Y}]
= 0$. Then $ G$ is metabelian.
\end{Lemma}
\begin{proof}
Since $\hat {X}$ and $\hat {Y}$ are projectively invariant by $ G$,
for each   $\hat g\in G$ there are constants $a_{\hat g}, b_{\hat
g}\in\mathbb{C}$ such that $\hat g_*\hat {X} = a_{\hat g}\cdot\hat
{X}$ and $\hat g_*\hat {Y} = b_{\hat g}\cdot \hat {Y}$. Given now a
tangent to the identity element $\hat h \in G_{\Id}$ we have $a_h=1$
and $b_h=1$,
 so that $\hat h_{*}\hat {X} = \hat {X}$ and  $\hat h_{*}\hat {Y} = \hat {Y}$.
  This implies that  $G_{\Id }$ is abelian. Since $DG$ is abelian,
  we have that $[G,G]< G_{\Id}$,
  so that $[G,G]$ is abelian.
\end{proof}

\begin{proof}[Proof of Theorem~\ref{Theorem:metabelian}]

Let $G < \hat{\Diff}(\mathbb{C}^{2},0)$ be a metabelian non-abelian
subgroup such that the group of commutators $[G,G]$ is  tangent to
the identity, $[G,G]< G_{\Id}$. If $G_{\Id}$ is trivial then the group is abelian. Therefore, we may assume that the Lie algebra of $G_{\Id}$ has positive dimension.  By
Proposition~\ref{proposition:abeliantangentidentity} we have two
cases:

\noindent{\bf Case 1}.  $[ G, G]\leq \langle \exp(N\hat {X})\mid N
{\rm \hspace{0.1cm}is \hspace{0.1cm} a
\hspace{0.1cm}rational\hspace{0.1cm} funtion}, \hat {X}(N) = 0
\rangle$ and $\hat f= \exp(\hat {X})\in [ G,  G]$. Then for all
$\hat g\in G$, $[\hat g, \hat f]\in[ G,  G]$ so that, there is a
rational function $\widetilde{N}$ such that $[\hat g, \hat f] =
\exp(\widetilde{N}\hat {X})$ then
 $\hat g \circ \exp(\hat {X})\circ \hat g^{-1} =
 \exp(\widetilde{N}\hat {X})\circ\exp(\hat {X}) =
\exp((\widetilde{N} + 1)\hat {X})$ therefore $\hat g^*(\hat {X}) =
N\hat {X}$ \vglue.1in

\noindent{\bf Case 2}.   $[ G, G]< \langle \exp(s\hat
{X})\circ\exp(t\hat {Y})\mid s,t\in\mathbb{C}\rangle$, take $\hat f=
\exp(\hat {X})$. Then for all $\hat g\in G$, $[\hat g, \hat f]\in[
G, G]$ so that, there are $\widetilde{s_1}$ and $t_1$ such that
$[\hat g, \hat f] = \exp(\widetilde{s_1}\hat {X})\circ\exp(t_1\hat
{Y})$ then $\hat g \circ \exp(\hat {X})\circ \hat g^{-1} =
\exp(\widetilde{s_1}\hat {X} + t_1\hat {Y})\circ \exp(\hat {X}) =
\exp(s_1\hat {X} + t_1\hat {Y})$. Therefore $\hat g^*(\hat {X}) =
s_1\hat {X} + t_1\hat {Y}$. Analogously we have $\hat g^*(\hat {Y})
= s_2\hat {X} + t_2\hat {Y}$.

This proves (ii) in Theorem~\ref{Theorem:metabelian}.

Let us now finish the proof, by
constructing in case (ii)  the formal closed meromorphic one-forms $\hat \omega_j, \,
j=1,2$. We can write  $\hat {X}_j = A_j\frac{\partial}{\partial x} +
B_j\frac{\partial}{\partial y}$ then
\[
\begin{bmatrix}
\frac{\partial g_1}{\partial x} & \frac{\partial g_1}{\partial y}  \\[7pt]
\frac{\partial g_2}{\partial x} & \frac{\partial g_2}{\partial y}
\end{bmatrix}
\begin{bmatrix}
A_1(z) & A_2(z) \\[7pt] B_1(z) & B_2(z)
\end{bmatrix}
=
\begin{bmatrix}
s_1A_1(\hat g) + t_1A_2(\hat g) & s_2A_1(\hat g) + t_2A_2(\hat g)\\[7pt]
s_1B_1(\hat g) + t_1B_2(\hat g) & s_2B_1(\hat g) + t_2B_2(\hat g)
\end{bmatrix}
\]
taking transposes, we have:
\[
\begin{bmatrix}
A_1(z) & B_1(z) \\[7pt] A_2(z) & B_2(z)
\end{bmatrix}
\begin{bmatrix}
\frac{\partial g_1}{\partial x} & \frac{\partial g_2}{\partial x}  \\[7pt]
\frac{\partial g_1}{\partial y} & \frac{\partial g_2}{\partial y}
\end{bmatrix}
=
\begin{bmatrix}
s_1A_1(\hat g) + t_1A_2(\hat g) & s_1B_1(\hat g) + t_1B_2(\hat g)\\[7pt]
 s_2A_1(\hat g) + t_2A_2(\hat g) & s_2B_1(\hat g) + t_2B_2(\hat g)
\end{bmatrix}
\]
thus
\[
\begin{bmatrix}
\frac{\partial g_1}{\partial x} & \frac{\partial g_2}{\partial x}  \\[7pt]
\frac{\partial g_1}{\partial y} & \frac{\partial g_2}{\partial y}
\end{bmatrix}
\begin{bmatrix}
s_1A_1(\hat g) + t_1A_2(\hat g) & s_1B_1(\hat g) + t_1B_2(\hat g)\\[7pt]
 s_2A_1(\hat g) + t_2A_2(\hat g) & s_2B_1(\hat g) + t_2B_2(\hat g)
\end{bmatrix}^{-1}
=
\begin{bmatrix}
A_1(z) & B_1(z) \\[7pt] A_2(z) & B_2(z)
\end{bmatrix}^{-1}
\]
so that we can  take:
\[
\begin{bmatrix}
C_1(z) & C_2(z) \\[7pt] D_1(z) & D_2(z)
\end{bmatrix}
=
\begin{bmatrix}
s_1A_1(z) + t_1A_2(z) & s_1B_1(z) + t_1B_2(z)\\[7pt]
s_2A_1(z) + t_2A_2(z) & s_2B_1(z) + t_2B_2(z)
\end{bmatrix}^{-1}
\]

and define  $\hat \omega_j = C_jdx + D_jdy$. Then $\hat g^*(\hat
\omega_1) = \hat g^*(\frac{1}{rQ(z)}.(s_2B_1 + t_2B_2, -(s_1B_1 +
t_1B_2)) = \frac{1}{Q(z)}(B_2,-A_2) = s_1\hat \omega_1 + s_2\hat
\omega_2$, where $Q(x) = C_1(z).D_2(z) - C_2(z).D_1(z)$ and $r =
s_1t_2 -s_2t_1$. Analogously $\hat g^*(\hat \omega_2) = t_1\hat
\omega_1 + t_2\hat \omega_2$, clearly $\hat \omega_1$ and $\hat
\omega_2$ are $\mathbb C$-linearly independent. It remains to show
that the $\hat \omega_j$ are closed forms, i,e. $\frac{\partial
D_j}{\partial x} - \frac{\partial C_j}{\partial y} = 0$. Since
$[\hat {X}_1,\hat {X}_2] = 0$ then
\begin{eqnarray*}
\frac{\partial A_2}{\partial x}A_1 +
\frac{\partial A_2} {\partial y}B_1 & =  &
\frac{\partial A_1}{\partial x}A_2 + \frac{\partial A_1}{\partial y}B_2\\
\frac{\partial B_2}{\partial x}A_1 +
\frac{\partial B_2}{\partial y}B_1 & = &
\frac{\partial B_1}{\partial x}A_2 + \frac{\partial B_1}{\partial y}B_2
\end{eqnarray*}
Thus ousting the value of $C_j$ and $D_j$ and using the above equations we can conclude.
\end{proof}

Next we study the possible normal forms of groups as in the
conclusion of Theorem~\ref{Theorem:metabelian}.

\begin{Remark}[groups leaving invariant a linear system of closed forms]
\label{Remark:normalformdouble}{\rm
Let $ G < {\hat \Diff}(\mathbb{C}^2,0)$ be a subgroup of formal
diffeomorphisms of two variables, that preserves
 the coordinate axes $(x = 0)$  and $(y = 0)$. Suppose that
 we have
\begin{equation}
\label{equation:invariance} \hat g^*(\hat \omega_j) = a_j\hat
\omega_1 + b_j\hat \omega_2, \, \,\forall \hat g \in  G.
\end{equation}
where $a_j, b_j\in\mathbb{C}^*$ and  $\hat \omega_j$ is a closed
formal meromorphic one-form. A diffeomorphism $\hat g \in G$ writes
$\hat g(x,y) = (xu,yv)$ where $u, v \in \hat {\mathcal O}_2$
 We have the following possibilities for
$\hat \omega_1, \hat \omega_2$ in suitable formal coordinates:

\begin{itemize}
\item[{\rm(1)}] ({\bf simple poles case}) If both forms have simple poles along
the coordinate axes we can write
 $\hat \omega_1 =
\alpha_1\frac{dx}{x} +\beta_1\frac{dy}{y}$ and $\hat \omega_2 =
\alpha_2\frac{dx}{x} +\beta_2\frac{dy}{y}$. From
equation~\eqref{equation:invariance} we get
\[
\alpha_1\frac{dx}{x} +\beta_1\frac{dy}{y}  + \alpha_1\frac{du}{u}
+\beta_1\frac{dv}{v} = \hat g^*(\hat \omega_j) = (a_1\alpha_1 +
b_1\alpha_2)\frac{dx}{x} + (a_1\beta_1 + b_1\beta_2)\frac{dy}{y}
\]
In matrix form we have:
\[
\begin{bmatrix}
\alpha_1 & \beta_1 \\[7pt] \alpha_2 & \beta_2
\end{bmatrix}
\begin{bmatrix}
\frac{dx}{x} \\[7pt]
\frac{dy}{y}
\end{bmatrix}
+
\begin{bmatrix}
\alpha_1 & \beta_1 \\[7pt] \alpha_2 & \beta_2
\end{bmatrix}
\begin{bmatrix}
\frac{du}{u} \\[7pt]
\frac{dv}{v}
\end{bmatrix}
=
\begin{bmatrix}
a_1 & b_1 \\[7pt] a_2 & b_2
\end{bmatrix}
\begin{bmatrix}
\alpha_1 & \beta_1 \\[7pt] \alpha_2 & \beta_2
\end{bmatrix}
\begin{bmatrix}
\frac{dx}{x} \\[7pt]
\frac{dy}{y}
\end{bmatrix}
\]

Comparing the poles we obtain:
\[
\begin{bmatrix}
a_1 & b_1 \\[7pt] a_2 & b_2
\end{bmatrix}
= \Id \hspace{0.4cm} {\rm and}  \hspace{0.4cm}
\begin{bmatrix}
\frac{du}{u} \\[7pt]
\frac{dv}{v}
\end{bmatrix}
= 0
\]
Thus  $\hat g(x,y) = (xu_0,yv_0)$ with
 $u_0$ and $v_0$ constant. In this case the group $G$ is linear.

\item[{\rm(2)}] ({\bf Pure polar  case})
Assume now that $\hat\omega_j$ has poles of order higher than one
and  no residues. We can write $\hat \omega_1 = d(\frac{1}{x^ny^m})$
and $\hat \omega_2 = d(\frac{1}{x^py^q})$. Given now a
diffeomorphism $\hat g(x,y) = (xu,yv)$  in G from
equation~\eqref{equation:invariance} we have:

\[
d(\frac{1}{x^ny^mu^nv^m}) = a_1d(\frac{1}{x^ny^m}) +
b_1d(\frac{1}{x^py^q}) \hspace{0.4cm} {\rm and } \hspace{0.4cm}
d(\frac{1}{x^py^qu^pv^q}) = a_2d(\frac{1}{x^ny^m}) +
b_2d(\frac{1}{x^py^q})
\]
Thus
\[
\frac{1}{x^ny^mu^nv^m} = \frac{a_1}{x^ny^m} + \frac{b_1}{x^py^q} +
k_1 \hspace{0.4cm}{\rm and }\hspace{0.4cm} \frac{1}{x^py^qu^pv^q} =
\frac{a_2}{x^ny^m} + \frac{b_2}{x^py^q}
\]
now
\[
\frac{1}{u^nv^m} = a_1 + b_1x^{n-p}y^{m - q} + k_1x^ny^m
\hspace{0.4cm}{\rm and }\hspace{0.4cm} \frac{1}{u^pv^q} = a_2 +
b_2x^{p - n}y^{q - m} + k_2x^py^q
\]
thus
\[
u = \frac{(a_2 +  b_2x^{p - n}y^{q - m} + k_2x^py^q)^{\frac{m}{D}}}
{(a_1 + b_1x^{n-p}y^{m - q} + k_1x^ny^m)^{\frac{q}{D}}}
\hspace{0.4cm}{\rm and }\hspace{0.4cm} u = \frac{(a_1 +
b_1x^{n-p}y^{m - q} + k_1x^ny^m)^{\frac{p}{D}}} {(a_2 +  b_2x^{p -
n}y^{q - m} + k_2x^py^q)^{\frac{n}{D}}}
\]
where $D = qn - pm$, therefore we have:
\[
\hat g(x,y) = \left(x.\frac{(a_2 +  b_2x^{p - n}y^{q - m} +
k_2x^py^q)^{\frac{m}{D}}} {(a_1 + b_1x^{n-p}y^{m - q} +
k_1x^ny^m)^{\frac{q}{D}}}\,,\, \,  y. \frac{(a_1 + b_1x^{n-p}y^{m -
q} + k_1x^ny^m)^{\frac{p}{D}}} {(a_2 +  b_2x^{p - n}y^{q - m} +
k_2x^py^q)^{\frac{n}{D}}}\right)
\]

\end{itemize}
Other mixed cases are studied in the same way. } \end{Remark}

The following example contradicts Corollary 4.2 in \cite{Brochero}.

\begin{Example}
\label{Example:solvnotmeta} {\rm An example of a tangent to the
identity group $G < \hat{\Diff}_{\Id}(\mathbb{C}^{2},0)$, which is
solvable but not metabelian is $G = < (\hat h(x), \hat a(x) + \hat
b(x)y); \hat h\in H >$, where $H <\hat{\Diff}_{\Id}(\mathbb{C},0)$
is any metabelian tangent to the identity subgroup, $\hat a(x)\in
\mathbb C[[x]]$ has order greater than $2$ and $\hat b(x)\in \mathbb
C[[x]]$ is  a unit, $\hat b(0) = 1$.}
\end{Example}


\section{Nilpotent groups and subalgebras of vector fields}
\label{Section:nilpotent}

We study nilpotent  subgroups of
$\hat{\Diff}_{\Id}(\mathbb{C}^{n},0)$,  i.e., nilpotent groups of
maps tangent to the identity. For  $n=1$ it is  known that the
concepts of solvable  and metabelian group are equivalent.
Nevertheless, this is not true for dimension $n \geq 2$: take $G<
\Diff(\mathbb C^2,0)$ as (induced by) the nilpotent linear group of
upper triangular matrices. This group has length $2$. In general, it
is known that for  a linear solvable group $G < \GL(n,\mathbb{C})$,
the solvable length is bounded by  the {\it Newman function}
$\rho(n)$, where $\rho(n)\leq 2n$, in particular, $\rho(2) = 4$ and
$\rho(3) = 5$ \cite{Newman}. In this section we prove
Theorem~\ref{Theorem:nilpotent}, i.e., that every nilpotent
subalgebra $\mathfrak{l}$ of $\hat{\mathcal{X}}(\mathbb{C}^n,0)$ has
length $l(\mathfrak{l})$ at most $n$.

 Note that $\hat{\mathcal{X}}(\mathbb{C}^{n},0)\otimes \hat{K}_{n}$
  is a vector space of dimension $n$ over $\hat{K}_{n}$,
  where $\hat{K}_{n}$ is the fraction field of $\hat{\mathcal{O}}_n$.
  Denote by $\mathcal{R}$ the center of $\mathfrak{l}$ and  by $\{X_1,\ldots,\mathcal{X}_k\}$
  a basis of $\mathcal{R}\otimes \hat{K}_{n}$.

\begin{Lemma}
\label{algebra nil} Let $\mathfrak{l}$ a nilpotent subalgebra of
$\hat{\mathcal{X}}(\mathbb{C}^{n},0)$ with dimension $m$ over
$\hat{K}_{n}$. Then there is an ordered basis
$\{\mathcal{X}_{1},\ldots, \mathcal{X}_m\}$ for  $\mathfrak{l}$ over
$\hat{K}_{n}$ in the following sense $\mathcal{Z}\in\mathfrak{l}$:
\begin{enumerate}
\item
\[
[\mathcal{Z},\mathcal{X}_l] =
\begin{cases}
0,&\mbox{if}\quad l = 1,\ldots, k,\\
\sum\limits_{j=1} ^{r<l}v_j\mathcal{X}_j, &\mbox{if}\quad l = k +
1,\ldots, m.
\end{cases}
\]
\item
If $\mathcal{Z} =
\sum\limits_{j=1}^{l}u_j\mathcal{X}_j\in\mathfrak{l}$ then
$\mathcal{X}_j(u_r) = 0$ with $(j = 1,\ldots, r)$ and $(r =
1,\ldots, l)$.
\end{enumerate}
\end{Lemma}
\begin{proof}
As $\mathfrak{l}$ a nilpotent subalgebra we have that $k\geq 1$.
Take $\mathcal{S} = (\sum\limits_{j=1}
^{k}\hat{K}_{n}\mathcal{X}_j)\cap \mathfrak{l}$,  we have that
$\mathcal{R}\subset\mathcal{S}$ and $\mathcal{S}$ is an abelian
subalgebra of $\mathfrak{l}$. Now for all
$\mathcal{Z}\in\mathfrak{l}$, we have $[\mathcal{Z},
\sum\limits_{j=1} ^{k}w_j\mathcal{X}_j] = \sum\limits_{j=1}
^{k}[\mathcal{Z}, w_j\mathcal{X}_j] = \sum\limits_{j=1}
^{k}(\mathcal{Z}(w_j)\mathcal{X}_j - w_j[\mathcal{Z},\mathcal{X}_j])
= \sum\limits_{j=1} ^{k}\mathcal{Z}(w_j)\mathcal{X}_j$,
($[\mathcal{Z},\mathcal{X}_j] = 0$, because
$\mathcal{X}_j\in\mathcal{R}$, $j = 1,\ldots, k$) then
$[\mathcal{Z}, \sum\limits_{j=1}
^{k}w_j\mathcal{X}_j]\in\mathcal{S}$,  thus $\mathcal{S}$ is an
ideal of $\mathfrak{l}$. Therefore $\mathfrak{l}/\mathcal{S}$ is a
nilpotent Lie algebra, thus $\mathcal{R}_1$ the center of
$\mathfrak{l}_1 = \mathfrak{l}/\mathcal{S}$ is not trivial, i.e.,
there are $\mathcal{X}_{k +
1},\ldots,\mathcal{X}_{p_1}\in\mathfrak{l}\setminus\mathcal{S}$ such
that $\overline{\mathcal{X}}_{k +
1},\ldots,\overline{\mathcal{X}}_{p_{1}}\in\mathcal{R}_1$ are the
generators the basis of $\mathcal{R}_1\otimes \hat{K}_{n}$. Clearly
$\{X_1,\ldots,\mathcal{X}_{p_1}\}$ are linearly independent in
$\mathfrak{l}\otimes \hat{K}_{n}$. Now as
$\overline{\mathcal{X}}_l\in\mathcal{R}_1$ $(l = k + 1,\ldots,
p_1)$, we have $\sum\limits_{j=1} ^{k}f_{l,j}\mathcal{X}_j =
[\mathcal{Z},\mathcal{X}_l + \sum\limits_{j=1}
^{k}w_{l,j}\mathcal{X}_j] = [\mathcal{Z},\mathcal{X}_l] +
[\mathcal{Z},\sum\limits_{j=1} ^{k}w_{l,j}\mathcal{X}_j] =
[\mathcal{Z},\mathcal{X}_l] +  \sum\limits_{j=1}
^{k}\mathcal{Z}(w_{l,j})\mathcal{X}_j$, then
$[\mathcal{Z},\mathcal{X}_l] = \sum\limits_{j=1} ^{k}( f_{l,j} -
\mathcal{Z}(w_{l,j}))\mathcal{X}_j = \sum\limits_{j=1}
^{k}v_{l,j}\mathcal{X}_j$, $(l = k + 1,\ldots, p_1)$. On the other
hand, if $\mathcal{Z} = \sum\limits_{j=1} ^{p_1}u_j\mathcal{X}_j +
\sum\limits_{j = p_1} ^{m}u_j\mathcal{Y}_j\in \mathfrak{l}$ for any
$\mathcal{Y}_j\in\mathfrak{l}$ completing the basis,  we have
$\mathcal{X}_j(u_r) = 0$ for $j = 1,\ldots, k$ and $r = 1,\ldots,
m$, because $\mathcal{X}_j\in\mathcal{R}$, $j = 1,\ldots, k$. Now
for $l = k + 1,\ldots, p_1$, we have $\sum\limits_{s=1}
^{k}v_{l,s}\mathcal{X}_s = [\sum\limits_{s=1}^{l}u_s\mathcal{X}_s,
\mathcal{X}_l] =
-\sum\limits_{s=1}^{k}\mathcal{X}_l(u_s)\mathcal{X}_s +
\sum\limits_{s= k + 1}^{l} (u_s[\mathcal{X}_s,\mathcal{X}_l] -
\mathcal{X}_l(u_s)\mathcal{X}_s) =
-\sum\limits_{s=1}^{k}\mathcal{X}_l(u_s)\mathcal{X}_s +
\sum\limits_{s= k + 1}^{l} u_s(\sum\limits_{j=1}
^{k}f_{r,j}\mathcal{X}_j) - \sum\limits_{s= k +
1}^{l}\mathcal{X}_l(u_s)\mathcal{X}_s  = \sum\limits_{j=1}
^{k}((\sum\limits_{s= k + 1}^{l} u_sf_{s,j}) -
\mathcal{X}_l(u_j))\mathcal{X}_j - \sum\limits_{s= k +
1}^{l}\mathcal{X}_l(u_s)\mathcal{X}_s$, then $\sum\limits_{j=1}
^{k}((\sum\limits_{s= k + 1}^{l} u_sf_{s,j}) - \mathcal{X}_l(u_j) -
v_{l,s})\mathcal{X}_j - \sum\limits_{s= k +
1}^{l}\mathcal{X}_l(u_s)\mathcal{X}_s = 0$, and as the
$\mathcal{X}_j$ are linearly independent, we have
$\mathcal{X}_l(u_r) = 0$ for $r = k + 1,\ldots, p_1$ and $l =
1,\ldots, p_1$, (note that $\mathcal{X}_l(u_r)$ for $r = 1,\ldots,
k$ can be nonzero ). If $p_1 = m$ we have nothing more to proof, in
the other case take $\mathcal{S}_1 = (\sum\limits_{j=1}
^{p_1}\hat{K}_{n}\overline{\mathcal{X}}_j)\cap \mathfrak{l}_1$ we
have that $\mathcal{R}_1\subset\mathcal{S}_1$ and $\mathcal{S}_1$ is
an abelian subalgebra of $\mathfrak{l}_1$ and analogously
$\mathcal{S}_1$ is an ideal of $\mathfrak{l}_1$. Therefore
$\mathfrak{l}_1/\mathcal{S}_1$ is a nilpotent Lie algebra, thus
$\mathcal{R}_2$ the center of $\mathfrak{l}_2 =
\mathfrak{l}_1/\mathcal{S}_1$ is not trivial, i.e., there are
$\overline{\mathcal{X}}_{p_1 +
1},\ldots,\overline{\mathcal{X}}_{p_2}\in\mathfrak{l}_1\setminus\mathcal{S}_1$
such that $\overline{\overline{\mathcal{X}}}_{p_1 +
1},\ldots,\overline{\overline{\mathcal{X}}}_{p_{2}}\in\mathcal{R}_2$
are the generators the basis of $\mathcal{R}_2\otimes \hat{K}_{n}$.
Clearly $\{X_1,\ldots,\mathcal{X}_{p_2}\}$ are Linearly Independent
in $\mathfrak{l}\otimes \hat{K}_{n}$, as
$\overline{\overline{\mathcal{X}}}_l\in\mathcal{R}_1$ $(l = p_1 +
1,\ldots, p_2)$, we have $\sum\limits_{j=1}
^{p_1}f_{l,j}\mathcal{X}_j = \sum\limits_{j=k + 1}
^{p_1}g_{l,j}.(\mathcal{X}_j + \sum\limits_{r=1}
^{k}h_{l,r}\mathcal{X}_r) = \sum\limits_{j=k + 1}
^{p_1}g_{l,j}\overline{\mathcal{X}}_j  = [\mathcal{Z},\mathcal{X}_l
+ \sum\limits_{j=1} ^{p_1}w_{l,j}\mathcal{X}_j] =
[\mathcal{Z},\mathcal{X}_l] + [\mathcal{Z},\sum\limits_{j=1}
^{p_1}w_{l,j}\mathcal{X}_j]$, then $[\mathcal{Z},\mathcal{X}_l] =
\sum\limits_{j=1} ^{p_1}v_{l,j}\mathcal{X}_j$, $(l = p_1 + 1,\ldots,
p_2).$ Similarly we get the second item, repeating this process a
finite number of times $(m < n)$ the lemma is proved.

\end{proof}

As consequence this lemma we have:
\begin{proof}[Proof of Theorem~\ref{Theorem:nilpotent}]
It is enough to prove that given a nilpotent subalgebra
$\mathfrak{l}$  of $\hat{\mathcal{X}}(\mathbb{C}^n,0)$ the $\mathcal
L$ has length $l(\mathfrak{l})$ at most $n$. If the dimension of
$\mathcal{R}\otimes \hat{K}_{n}$ is $n$, then for all
$\mathcal{Z}\in\mathfrak{l}$, $\mathcal{Z} = \sum\limits_{j=1}
^{n}u_j\mathcal{X}_j$ and $0 = [\mathcal{Z},\mathcal{X}_k] =
[\sum\limits_{j=1} ^{n}u_j\mathcal{X}_j,\mathcal{X}_k] =
\sum\limits_{j=1} ^{n}[u_j\mathcal{X}_j,\mathcal{X}_k] =
\sum\limits_{j=1} ^{n}(u_j[\mathcal{X}_j,\mathcal{X}_k] -
\mathcal{X}_k(u_j)\mathcal{X}_j) = -\sum\limits_{j=1} ^{n}
\mathcal{X}_k(u_j)\mathcal{X}_j$, then $\mathcal{X}_k(u_j) = 0$ for
$(j,k = 1,\ldots, n)$. thus $u_j$ are constants and therefore
$\mathfrak{l}$ is an abelian Lie algebra. Now if the dimension of
$\mathcal{R}\otimes \hat{K}_{n}$ is $m$, where $m$ is the dimension
of $\mathfrak{l}\otimes \hat{K}_{n}$ we have nothing to proof.
Finally  if the dimension of $\mathcal{R}\otimes \hat{K}_{n}$ is $k
< m$ by the Lemma \ref{algebra nil} we have for
$\mathcal{Z}_1,\mathcal{Z}_2\in \mathfrak{l}$,
$[\mathcal{Z}_1,\mathcal{Z}_2] =
[\sum\limits_{j=1}^{m}u_j\mathcal{X}_j,\sum\limits_{r=1}
^{m}v_r\mathcal{X}_r] = \sum\limits_{j=1} ^{m}\sum\limits_{r=1}
^{m}[u_j\mathcal{X}_j,v_r\mathcal{X}_r] = \sum\limits_{j=1}
^{m}\sum\limits_{r=1} ^{m}(u_j\mathcal{X}_j(v_r)\mathcal{X}_r -
v_r\mathcal{X}_r(u_j)\mathcal{X}_j +
u_j.v_r.[\mathcal{X}_j,\mathcal{X}_r]) = \sum\limits_{j=1} ^{m -
1}w_j\mathcal{X}_j$, because $\mathcal{X}_j(u_{m}) = 0$ $(j =
1,\ldots, m)$ and $[\mathcal{Z},\mathcal{X}_l] = \sum\limits_{j=1}
^{l - 1}v_j\mathcal{X}_j$. Thus all $\mathcal{Z}\in\mathfrak{l}^1$
is the form $\mathcal{Z} = \sum\limits_{j=1} ^{m -
1}w_j\mathcal{X}_j$. Now if $\mathcal{Z}_1,\mathcal{Z}_2\in
\mathfrak{l}^1$ we have $[\mathcal{Z}_1,\mathcal{Z}_2] =
[\sum\limits_{j=1} ^{m - 1}u_j\mathcal{X}_j,\sum\limits_{r=1} ^{m -
1}v_r\mathcal{X}_r] = \sum\limits_{j=1} ^{m - 1}\sum\limits_{r=1}
^{m - 1}[u_j\mathcal{X}_j,v_r\mathcal{X}_r] = \sum\limits_{j=1} ^{m
- 1}\sum\limits_{r=1} ^{m - 1}(u_j\mathcal{X}_j(v_r)\mathcal{X}_r -
v_r\mathcal{X}_r(u_j)\mathcal{X}_j +
u_j.v_r.[\mathcal{X}_j,\mathcal{X}_r]) = \sum\limits_{j=1} ^{m -
2}w_j\mathcal{X}_j$, by the Lemma \ref{algebra nil}. Then all
$\mathcal{Z}\in\mathfrak{l}^2$ is the form $\mathcal{Z} =
\sum\limits_{j=1} ^{m - 2}w_j\mathcal{X}_j$. Repeating this process
at most $m - 2$ times we can conclude.
\end{proof}

\begin{Remark}
{\rm In \cite{Martelo-Ribon} a detailed study of the length  for
solvable subgroups of $\hat{\Diff}(\mathbb C^n,0)$ is found.  The
authors correct a statement from \cite{Brochero} and prove the
following more general result: {\it Let $G < \hat{\Diff}(\mathbb C^n,0)$ be a solvable
group. Then the soluble length of  $G$ is at most $2n -1 + \rho(n)$
where $\rho : {\mathbb N} \to {\mathbb N}$ is the Newman function.}

}

\end{Remark}


\vglue.2in
\centerline{\textsc{ Part II - Groups with dicritical maps}}
\section{Dicritic groups with abelian commutators}

Unlike the one-dimensional case two commuting tangent to the
identity diffeomorphisms may have different orders of tangency to
the identity:  take  $\hat f = \exp(\hat {X})$ and $\hat g =
\exp(\hat {Y})$, where the vector fields $\hat {X}$ and $\hat {Y}$
are  given as in $(1)$ above. This is the main reason why we do not
have an equivalence between the concepts of metabelian,
quasi-abelian and solvable groups in dimension $n \geq 2$. From now
on we shall   take a closer look at this issue. Firstly, in this
section, we investigate the characterization of quasi-abelian
groups. For this we shall refer to the following concepts, which are
two main notions in this paper.

Now we pave the way to Theorem~\ref{Theorem:C}. For the first part
we shall need some lemmas below.
\begin{Lemma}
\label{Radial} Let $\hat {X} = f(z)\vec{R}$ and $\hat {Y} =
g(z)\vec{R}$, where $\vec R$ is the radial vector field and $f$ and
$g$ are homogeneous polynomials of degree $k$ e $s$ respectively.
Then $[\hat {X},\hat {Y}] = 0$ if and only if $k = s$.
\end{Lemma}
\begin{proof}
Trivial, because $[\hat {X},\hat {Y}] = (k - s)f(z)g(z)\vec{R}$.
\end{proof}

\begin{Lemma}
\label{dicritico, qualq ----- dicrit} Let $\hat {X}\in\hat{\mathcal
X}_{k + 1}(\mathbb{C}^n,0)$  a dicritic vector field. For any vector
field $\hat {Y}$ with order greater than $2$, such that $[\hat
{X},\hat {Y}] = 0$, we have that $\hat {Y}$ is a dicritic vector
field with order $k + 1$.
\end{Lemma}
\begin{proof}
Suppose that $\hat {Y}$ has order $r\geq 2$, thus:
\[
\hat {X} = f(z)\overrightarrow{R} + (p_{k + 2}^{(1)} +
\cdots)\frac{\partial}{\partial z_{1}} + \cdots + (p_{k + 2}^{(n)} +
\cdots)\frac{\partial}{\partial z_{n}}
\]
\[
\hat {Y} = (q_{r}^{(1)} + \cdots)\frac{\partial}{\partial z_{1}} +
\cdots + (q_{r}^{(n)} + \cdots)\frac{\partial}{\partial z_{n}}
\]

Now, the term of lower order of $[\hat {X},\hat {Y}]$ is:
\begin{align*}
[f\overrightarrow{R},q_{r}^{(1)}\frac{\partial}{\partial z_{1}}
 + \cdots + q_{r}^{(n)}\frac{\partial}{\partial z_{n}}] =
 (r - 1)f.(q_{r}^{(1)}\frac{\partial}{\partial z_{1}}
 &+ \cdots + q_{r}^{(n)}\frac{\partial}{\partial z_{n}}) -
 (q_{r}^{(1)}\frac{\partial f}{\partial z_{1}}
 +\\
 &\cdots + q_{r}^{(n)}\frac{\partial f}{\partial z_{n}}).\overrightarrow{R}
 \end{align*}
As $[\hat {X},\hat {Y}] = 0$, we have that $(r - 1)f.q_{r}^{(j)} =
(\nabla f.Q_{r})z_j$, onde $Q_{r} = (q_{r}^{(1)}, \ldots,
q_{r}^{(n)})$. Thus $(r - 1)f.q_{r}^{(1)}z_j = (r -
1)f.q_{r}^{(j)}z_1$ and as $q_{r}^{(1)}\neq 0$, we have $q_{r}^{(j)}
= \frac{q_{r}^{(1)}}{z_1}.z_j = g.z_j$, for $j= 1,\ldots,n$. So the
1-Jet of $\hat {Y}$ is $g\overrightarrow{R}$, therefore $\hat {Y}$
is dicritic vector  field and by the previous lemma $\hat {Y}$ have
order $k + 1$.
\end{proof}

\begin{Lemma}
\label{campos multiplos} Let $\hat {X}, \hat {Y} \in \hat{\mathcal
X}_{k}(\mathbb{C}^n,0)$, $k\geq 2$. Suppose that $\hat {X}$ is
regular dicritic and $\hat {Y}$ is dicritic. If $[\hat {X},\hat {Y}]
= 0$,  then there is $c\in\mathbb{C}\setminus\{0\}$ such that $\hat
{Y} = c.\hat {X}$.
\end{Lemma}
\begin{proof}
Since $\hat {X}$ and $\hat {Y}$ are dicritic, then
\[
\hat {X} = f(z)\vec{R} + (p_{k + 2}^{(1)} +
\cdots)\frac{\partial}{\partial z_{1}} + \cdots + (p_{k + 2}^{(n)} +
\cdots)\frac{\partial}{\partial z_{n}}
\]
\[
\hat {Y} = g(z)\vec{R} + (q_{k + 2}^{(1)} +
\cdots)\frac{\partial}{\partial z_{1}} + \cdots + (q_{k + 2}^{(n)} +
\cdots)\frac{\partial}{\partial z_{n}}
\]
We have $[f\vec{R},g\vec{R}] = 0$, by Lemma \ref{Radial}. Since
$[\hat {X},\hat {Y}] = 0$, then the $2k + 2$-jet to Lie bracket is
\[
[g\vec{R},p_{k + 2}^{(1)}\frac{\partial}{\partial z_{1}}
 + \cdots + p_{k + 2}^{(n)}\frac{\partial}{\partial z_{n}}] -
 [f\vec{R},q_{k + 2}^{(1)}\frac{\partial}{\partial z_{1}}
 + \cdots + q_{k + 2}^{(n)}\frac{\partial}{\partial z_{n}}] = 0
\]
Now note that
\begin{align*}
[f\vec{R},q_{k + 2}^{(1)}\frac{\partial}{\partial z_{1}}
 + \cdots + q_{k + 2}^{(n)}\frac{\partial}{\partial z_{n}}] =
 (k + 1)f.(q_{k + 2}^{(1)}\frac{\partial}{\partial z_{1}}
 &+ \cdots + q_{k + 2}^{(n)}\frac{\partial}{\partial z_{n}}) -
 (q_{k + 2}^{(1)}\frac{\partial f}{\partial z_{1}}
 +\\
 &\cdots + q_{k + 2}^{(n)}\frac{\partial f}{\partial z_{n}}).\vec{R}
 \end{align*}
Then we have
\[
(k + 1)(f.q_{k + 2}^{i} - g.p_{k + 2}^{i}) = z_{i}(\nabla f.Q_{k +
2} - \nabla g.P_{k + 2})
\]
for $i\in\{1,\ldots,n\}$, thus
\[
\frac{f.q_{k + 2}^{i_{0}} - g.p_{k + 2}^{i_{0}}}{z_{i_{0}}} =
\frac{f.q_{k + 2}^{j_{0}} - g.p_{k + 2}^{j_{0}}}{z_{j_{0}}}
\]
or equivalently, $f.(z_{j_{0}}q_{k + 2}^{(i_{0})} - z_{i_{0}}q_{k +
2}^{(j_{0})}) = g.(z_{j_{0}}p_{k + 2}^{(i_{0})} - z_{i_{0}}p_{k +
2}^{(j_{0})})$. But by hypothesis $f$ and $z_{j_{0}}p_{k +
2}^{(i_{0})} - z_{i_{0}}p_{k + 2}^{(j_{0})}$ are coprime, then $f\mid g$.
As $f$ and $g$ has the same degree $g = c.f$ were
$c\in\mathbb{C}^{*}$. Thus the $2k + 2$- jet of Lie bracket is:
\[
[f\overrightarrow{R},(q_{k + 2}^{(1)} - cp_{k +
2}^{(1)})\frac{\partial}{\partial z_{1}}
 + \cdots + (q_{k + 2}^{(n)} - cp_{k + 2}^{(n)} )\frac{\partial}{\partial z_{n}}] = 0
\]
Using the same argument of the previous lemma we have
\[
(q_{k + 2}^{(1)} - cp_{k + 2}^{(1)})z_j = (q_{k + 2}^{(j)} - cp_{k +
2}^{(j)})z_1
\]
so,  $q_{k + 2}^{(1)} - cp_{k + 2}^{(1)} = 0$, in consequence $q_{k
+ 2}^{(j)} - cp_{k + 2}^{(j)} = 0$, for all $j = 1,\ldots,n$, or
\[
(q_{k + 2}^{(1)} - cp_{k + 2}^{(1)})\frac{\partial}{\partial z_{1}}
 + \cdots + (q_{k + 2}^{(n)} - cp_{k + 2}^{(n)} )\frac{\partial}{\partial z_{n}} =
\frac{(q_{k + 2}^{(1)} - cp_{k + 2}^{(1)})}{z_1}\overrightarrow{R}
\]
but this latter does not occur, because by the Lemma \ref{Radial} we
have $\frac{(q_{k + 2}^{(1)} - cp_{k + 2}^{(1)})}{z_1}$ has degree
$k$, and this is impossible because $q_{k + 2}^{(1)} - cp_{k +
2}^{(1)}$ has degree $k + 2$.
Then, we have $q_{k + 2}^{(j)} = cp_{k + 2}^{(j)}$, $\forall j\in \{ 1,\ldots,n\}$.\\
Finally suppose that $Q_{k + j} = cP_{k + j}$ for $j = 1,\ldots,i$,
the $(2k + i + 1)$-jet of Lie bracket is
\[
[g\vec{R},p_{k + i + 1}^{(1)}\frac{\partial}{\partial z_{1}} +
\cdots + p_{k + i + 1}^{(n)}\frac{\partial}{\partial z_{n}}] -
[f\vec{R},q_{k + i + 1}^{(1)}\frac{\partial}{\partial z_{1}} +
\cdots + q_{k + i + 1}^{(n)}\frac{\partial}{\partial z_{n}}] = 0
\]
by the supposed the following sum is symmetric
\[
\sum\limits_{j = 2}^{i} [p_{k + j}^{(1)}\frac{\partial}{\partial
z_{1}} + \cdots + p_{k + j}^{(n)}\frac{\partial}{\partial z_{n}},
q_{k + i + 2 - j}^{(1)}\frac{\partial}{\partial z_{1}} + \cdots +
q_{k + i + 2 - j}^{(n)}\frac{\partial}{\partial z_{n}}] = 0
\]
Then similarly to the case $k + 2$ we have that $Q_{k + j + 1} =
cP_{k + j + 1}$ therefore $\hat {Y} = c\hat {X}$
\end{proof}

\begin{Remark}
{\rm We cannot exclude  the regularity condition in the previous
lemma, since the two vector fields of item $2.(d)$ in
Example~\ref{Example:abelianformalinvariant}  are dicritic and
commute, but they are not regular dicritic and are not $\mathbb
C$-linearly dependent.}
\end{Remark}

The following proposition is found in \cite{Brochero}.
\begin{Proposition}[\cite{Brochero} Proposition 4.2]
\label{dicfabio} Let $\hat f\in \hat{\Diff}_{r +
1}(\mathbb{C}^{n},0)$ and $\hat g \in \hat{\Diff}_{s +
1}(\mathbb{C}^{n},0)$. Suppose that $\hat f$ is dicritic and
commutes with $\hat g$. Then $r = s$ and $G$ is also dicritic.
\end{Proposition}

We state now  the main tool in the proof of Theorem~\ref{Theorem:C}.

\begin{Proposition}
\label{comu dicr} Let $ G < \hat{\Diff}_{\Id}(\mathbb{C}^{n},0)$ be
a subgroup of diffeomorphisms tangent to the identity and  $\hat
f\in \hat{\Diff}_{\Id}(\mathbb{C}^{n},0)$ a regular dicritic
diffeomorphism. If $\hat f$ commutes with $G$ then $ G <
\langle \exp(t\hat X)\mid t\in\mathbb{C}\rangle$, where $\hat f =
\exp(\hat {X})$. In particular, $ G $ is abelian.
\end{Proposition}
\begin{proof}
Let $\hat g \in  G $ be a diffeomorphism, from Proposition
\ref{dicfabio}, $\hat g$ is a dicritic diffeomorphism of  same order
than $\hat f$, we say $k + 1$. From  the exponential bijection there
is $\hat {Y}$ such that $\exp(\hat {Y})= \hat g$.  Then
\[
\hat {Y} = g(z)\vec{R} + (q_{k + 2}^{(1)} +
\cdots)\frac{\partial}{\partial z_{1}} + \cdots + (q_{k + 2}^{(n)} +
\cdots)\frac{\partial}{\partial z_{n}}
\]
and thus $\hat {Y}$ is dicritic. Since $\hat f$ commutes with $
\hat g $ it commutes with $\hat g(z) = \exp(\hat {Y})(z)$ and from
Lemma \ref{fluxo}, $\hat f$ commutes with $\exp(t\hat {Y})(z)$ for
all $t\in\mathbb{C}$. Similarly, for each $t$  we have that
$\exp(t\hat {Y})(z)$ commutes with $\exp(s\hat {X})(z)$ and thus
$[\hat {X},\hat {Y}] = 0$. From Lemma \ref{campos multiplos},
there is $r\in\mathbb{C}$ such that $\hat {Y} = r\hat {X}$.
Consequently $\hat g(z) = \exp(r\hat {X})(z)$ and therefore $ G <
\langle \exp(t\hat {X})\mid t\in\mathbb{C}\rangle$.
\end{proof}

 The above proposition motivates the following:
\begin{Definition}[Unimodular diffeomorphism]
{\rm A formal diffeomorphism tangent to the identity $\hat f\in \hat{\Diff}(\mathbb C^n,0)$
is {\it unimodular} if its centralizer in $\hat{\Diff}_{\Id}$ is $\mathbb C=\{\hat f^{[t]},
\, t \in \mathbb C\}$, i.e., generated by $\hat f$. }
\end{Definition}

Proposition~\ref{comu dicr} then says that a regular dicritic formal diffeomorphism is unimodular.

\begin{proof}[Proof of Theorem~\ref{Theorem:C}]
Let $ G < {\hat \Diff}(\mathbb{C}^{n},0)$ be a subgroup containing  a
regular dicritic diffeomorphism $\hat f=\exp(\hat {X})\in G $.
First, suppose that $ G $ is quasi-abelian.  From Proposition
\ref{comu dicr}, we have $ G _{\Id}< \langle \exp(t\hat {X})\mid
t\in\mathbb{C}\rangle$. Let $\hat g \in G $, as $\hat f\in G_{\Id}$
then $[\hat f,\hat g]\in  G_{\Id}$. Thus there is $t_{\hat
g}\in\mathbb{C}^{*}$ such that $[\hat f, \hat g] = \exp(t_{\hat
g}\hat X)$, then $\hat g \circ \hat f\circ \hat g ^{-1}\circ \hat f^{-1}
= \exp(t_{\hat g}\hat {X})$ so that
\begin{align*}
\hat g \circ \exp(\hat {X})\circ \hat g ^{-1} &= \exp(t_{\hat g}\hat
X)\circ \exp(\hat {X}) =
\exp((t_{\hat g} + 1)\hat {X})\\
&= \exp(c_{\hat g}\hat {X})
\end{align*}
from the same argument used in the proof of Lemma \ref{fluxo}, we
have $\forall s\in\mathbb{C}$, $\hat g \circ \exp(s\hat {X})\circ
\hat g ^{-1} = \exp(sc_{\hat g}\hat {X})$. Therefore $\hat g^{*}\hat
X = c_{\hat g}\hat {X}$, $\forall \hat g \in G $. Conversely,
suppose that $\hat X$ is projective invariant with respect to the
group $G$.  Therefore, for each   $ \hat g \in G $ there exists  $
c_{\hat g}$ such that $\hat g^{*}\hat {X} = c_{\hat g}\hat {X}$. We
claim that $\forall \hat g \in G_{\Id}$, $c_{\hat g} = 1$. In fact,
if $\hat f(z) = z + f(z)z + \cdots$ then $\exp (c_{\hat g}\hat
{X})(z) = z + c_{\hat g}f(z)z + \cdots$. Thus if $\hat g \in
G_{\Id}$, $\hat g \circ \hat f(z) = z + f(z)z + \hat g_{k+1}(z) +
\cdots$ and $\exp (c_{\hat g}\hat {X})(z)\circ \hat g(z) = z +
c_{\hat g}.f(z)z + \hat g_{k+1}(z) + \cdots$, then $c_{\hat g} = 1$.
Consequently $\forall \hat g \in G_{\Id}$, $\hat g^{*}\hat {X} =
\hat {X}$, {\it i.e.}, $ G_{\Id}$ commutes with $\hat f$. From
Proposition \ref{comu dicr}, $ G_{\Id}$ is abelian, {\it i.e.}, $ G
$ is quasi-abelian.
\end{proof}

In the same way as Theorem~\ref{Theorem:C} we have:

\begin{proof}[Proof of Proposition~\ref{Proposition:exttheorem2}]
It is immediate to verify that $(1) \Rightarrow (2)$. Let us now
prove $(2) \Rightarrow (1)$. Since $D G < G $, for all $\hat g \in G
$ we have that $\tilde g = D\hat g ^{-1}(0)\circ \hat g \in G_{\Id}$.
From $(2)$ we have that $ G $ commutes with $\hat f$ then $D\hat g
^{-1}(0)\circ \hat g = \exp(c_{\tilde g}X)$, therefore $G = D\hat
g(0)\circ\exp(c_{\tilde g}X)$, $\forall \hat g \in G $. Now let
$\hat g, \hat h \in G $ be diffeomorphisms, as $D G < G $ and from
$(2)$, we have that:
\[
\hat g \circ \hat h = D\hat g(0)\circ\exp(c_{\tilde g}X)\circ D\hat
h(0)\circ\exp(c_{\tilde h}X) = \hat h\circ \hat g
\]
Therefore $ G $ is abelian.
\end{proof}

\section{Metabelian  groups with dicritic elements}

Now we shall study metabelian   subgroups of $\hat\Diff(\mathbb C^n,
0)$ and prove Theorem~\ref{Theorem:metabeliandicritic}. We strongly
rely  on the preceding argumentation. The main step is:

\begin{Proposition}
\label{proposition:theorem3}
 Let $ G < \hat{\Diff}(\mathbb{C}^{n},0)$ be a subgroup with
 $D G $ abelian, and  $\hat f=\exp(\hat {X})\in G $  a regular
 dicritic diffeomorphism. Then $\hat {X}$ is projectively
 invariant by $ G $ if, and only if,  $\hat f$ commutes with $[G,G]$.
\end{Proposition}

\begin{proof}
Since  $D G $ is abelian the subgroup of commutators     of $ G $ is tangent
to the identity {\it i,e.}, $[ G , G ]< G_{\Id}$. If  $\hat f$
commutes with $[ G , G ]$, from
Proposition~\ref{comu dicr} we have
that $[ G , G ]< \langle \exp(t\hat {X})/t\in\mathbb{C}\rangle$
and from the proof of Theorem~\ref{Theorem:C} we know that $\hat {X}$ is
projectively invariant by $ G $. Assume now that  the vector field
$\hat {X}$ is projectively invariant by $ G $. Once again, since $[
G ,  G ]<  G_{\Id}$ we know that if $S\in[ G , G ]$ then $c_{S} =
1$. Thus $\forall S\in[ G , G ]$, $S^{*}\hat {X} = \hat {X}$,
therefore $\hat f= \exp(\hat {X})$ commutes with $[ G , G ]$.
\end{proof}

\begin{proof}[Proof of Theorem~\ref{Theorem:metabeliandicritic}]
Notice that if $G$ is quasi-abelian and  $DG$ is abelian then $G$ is
metabelian. Therefore, Theorem~\ref{Theorem:metabeliandicritic}
follows from Propositions~\ref{comu dicr} and
~\ref{proposition:theorem3}.
\end{proof}

Next we prove the announced  partial converse of $(2)$ in Theorem~\ref{Theorem:metabeliandicritic}.
\begin{proof}[Proof of Proposition~\ref{Proposition:converseTheoremB}]
 By hypothesis  $[ G , G ]$ is an abelian subgroup of
  diffeomorphisms tangent to the identity and we have that
\[
[\hat f,\hat h](z) = z + \lambda(\lambda^{k} - 1)f(z)z +
\lambda(\lambda^{k + 1} - 1)P_{k + 2}(z) + \cdots
\]
thus $[\hat f,\hat h] = \exp(\hat {X})$, where
\[
\hat {X} = \lambda(\lambda^{k} - 1)f(z)\vec{R} + (\lambda(\lambda^{k
+ 1} - 1)p_{k + 2}^{(1)} + \cdots)\frac{\partial}{\partial z_{1}} +
\cdots + (\lambda(\lambda^{k + 1} - 1)p_{k + 2}^{(n)} +
\cdots)\frac{\partial}{\partial z_{n}}
\]
since $\hat f$ is regular dicritic and $\lambda^k \ne 1,
\lambda^{k+1} \ne 1$, the above expression implies that $[\hat
f,\hat h]$ is regular dicritic. According to  Proposition
\ref{proposition:theorem3}  there is a projectively invariant formal
vector field $\hat {X} = \exp([\hat f,\hat h])$.
\end{proof}

Now we give an application of our results:

\begin{Corollary}
\label{corollary:theoremB} Let $ G  = \langle \hat f, \hat h\rangle
< {\hat \Diff}(\mathbb{C}^{n},0)$, where $\hat f$ is regular
dicritic and $\hat h(z) = \lambda z$, $\lambda^{k} \neq 1$,
$\lambda^{k + 1} \neq 1$. The group $ G $ is metabelian if and only
if $[\hat f,\hat h^{2}]$ and $[\hat f^{2},\hat h]$ commute with
$[\hat f,\hat h]$.
\end{Corollary}
\begin{proof}
Let us first introduce some notation.
Given two elements $\hat \varphi, \hat \psi \in \hat \Diff(\mathbb
C^2,)$ we shall write $\hat \varphi^* \hat \psi := \hat \varphi
\circ \hat \psi \circ \hat \varphi^{-1}$.

If $ G < \hat \Diff(\mathbb
C^2,0)$ is metabelian,  it is immediately  seen that $[\hat f,\hat
h^{2}]$ and  $[\hat f^{2},\hat h]$ commute with $[\hat f,\hat h]$.
We prove the converse, in fact we have that $[\hat f,\hat h]$ is
regular dicritic as noted above. Now, as $[\hat f,\hat h^{2}]$,
$[\hat f^{2},\hat h]$ commutes with $[\hat f,\hat h]$, we have
$[\hat f^{2},\hat h] = \exp(c\hat {X})$ and $[\hat f,\hat h^{2}] =
\exp(r\hat {X})$, since $\hat f^{*}[\hat f,\hat h]\circ [\hat f,\hat
h] = [\hat f^{2},\hat h]$ and $[\hat f,\hat h]\circ \hat h^{*}[\hat
f,\hat h] = [\hat f,\hat h^{2}]$ then
\[
\exp(c\hat {X}) = [\hat f^{2},\hat h] = \hat f^{*}[\hat f,\hat
h]\circ [\hat f,\hat h] = \hat f^{*}\exp(\hat {X})\circ \exp(\hat {X})
\]
consequently, $\hat f^{*}\exp(\hat {X}) =  \exp(c\hat {X})\circ
\exp(-\hat{X}) = \exp(\widetilde{c}\hat {X})$. Using the same
argument in the proof of Lemma \ref{fluxo}, $\hat f^{*}\hat {X} =
\widetilde{c}\hat {X}$. Similarly $\hat h^{*}\hat {X} =
\widetilde{r}\hat {X}$. Thus by
Theorem~\ref{Theorem:metabeliandicritic} the group  $ G $ is
metabelian.
\end{proof}

\section{Solvable groups with some dicritic element}

The next three lemmas will be used in the proof of  Theorem~\ref{Theorem:D}.

\begin{Lemma}
\label{comutador} Let $\hat f \in {\hat \Diff}_{r}(\mathbb{C}^{n},0)$ and $\hat g \in  {\hat \Diff}_{s}(\mathbb{C}^{n},0)$ be formal diffeomorphisms. Then
\[
\hat f(\hat g(z))- \hat g (\hat f(z))= D\hat f_{r + 1}(z)\hat g_{s +
1}(z) - D\hat g_{s + 1}\hat f_{r + 1}(z) + O(|z|^{r + s + 2})
\]
so $[\hat f, \hat g]= \Id$ or $[\hat f, \hat g] \in {\hat
\Diff}_{p}(\mathbb{C}^{n},0)$ with $p \geq r + s$.
\end{Lemma}
\begin{proof}
Let $\hat f(z)= z +  \sum\limits_{k=r} ^{r+s}\hat f_{k + 1}(z)+
O(|z|^{r + s + 2})$ and
 $\hat g(z)= z + \sum\limits_{j=s} ^{r+s}\hat g_{j + 1}(z)+  O(|z|^{r + s + 2})$ then:
\begin{align*}
 \hat f(\hat g(z)) &= z + \sum\limits_{j=s}^{r+s}
 \hat g_{j+ 1}(z)+ O(|z|^{r + s + 2}) +\\
 &\hspace{0.8cm}+ \sum\limits_{k=r}^{r+s}\hat f_{k + 1}
 (z + \sum\limits_{j=s}^{r+s}\hat g_{j + 1}(z) +
 O(|z|^{r + s + 2})) + O(|z|^{r + s + 2})  \\
&= z + \sum\limits_{j=s} ^{r+s}\hat g_{j+ 1}(z)+ \sum\limits_{k=r}
^{r+s}(\hat f_{k+ 1}(z)+
D\hat f_{k + 1}\hat g_{s + 1}(z) + O(|z|^{k + s + 2})) \\
&\hspace{2.1cm} +  O(|z|^{r + s + 2}) \\
&= z + \sum\limits_{j=s} ^{r+s}\hat g_{j+ 1}(z) + \sum\limits_{k=r}
^{r+s}\hat f_{k+ 1}(z) + D\hat f_{r + 1}(z)\hat g_{s + 1}(z) +
O(|z|^{r + s + 2})
\end{align*}
Similarly we have:

$
\hat g (\hat f(z))= z +\sum\limits_{k=r} ^{r+s}\hat f_{k+ 1}(z)
 + \sum\limits_{j=s} ^{r+s}\hat g_{j+ 1}(z) +
 D\hat g_{s + 1}(z)\hat f_{r + 1}(z) +  O(|z|^{r + s + 2}) \\
$
subtracting these two equalities we get the lemma.
\end{proof}

\begin{Lemma}
Let $\hat f \in \hat{\Diff}_{r}(\mathbb{C}^{n},0)$ and $\hat g
\in \hat{\Diff}_{s}(\mathbb{C}^{n},0)$ be dicritic
diffeomorphisms with $r\neq s$, given by
\[
\hat f(z) = z + f(z)z +  \cdots \hspace{1.0cm}\hat g(z) = z + g(z)z
+ \cdots
\]
then $[\hat f,\hat g]\in  \hat{\Diff}_{s + r}(\mathbb{C}^{n},0)$
is dicritic and given by
\[
[\hat f,\hat g](z)  = z + (r - s)g(z)f(z)z + \cdots
\]
\end{Lemma}
\begin{proof}
From Lemma \ref{comutador} the term of smaller order of $[\hat
f,\hat g]$ is
\begin{align*}
D\hat f_{r + 1}(z)\hat g_{s + 1}(z) - D\hat g_{s + 1}\hat f_{r +
1}(z) &= ( f(z)I + (z_{i}\frac{\partial f}{\partial z_{j}}))\hat
g_{s + 1}(z) -
D\hat g_{s + 1}f(z)z\\
&= ( f(z)I + (z_{i}\frac{\partial f}{\partial z_{j}}))
\hat g_{s + 1}(z) - (s + 1)f(z)\hat g_{s + 1}\\
&= ( -sf(z)I + (z_{i}\frac{\partial f}{\partial z_{j}}))\hat g_{s + 1}(z)\\
&= ( -sf(z)I + (z_{i}\frac{\partial f}{\partial z_{j}}))g(z)z
\end{align*}
Then the i-th component of this is
\[
-sf(z)g(z)z_{i} + g(z)z_{i}\nabla f(z) z = -sf(z)g(z)z_{i} +
rf(z)g(z)z_{i}
\]
thus
\[
D\hat f_{r + 1}(z)\hat g_{s + 1}(z) - D\hat g_{s + 1}\hat f_{r +
1}(z) = (r - s)f(z)g(z)z
 \]
therefore $[\hat f,\hat g]\in  \hat{\Diff}_{s + r +
1}(\mathbb{C}^{n},0)$ and this is dicritic, given by
\[
\hat h(z) = z + (r - s)g(z)f(z)z + \cdots
\]
\end{proof}

\begin{Lemma}
If a subgroup  $ G < \hat{\Diff}(\mathbb{C}^{n},0)$ contains two elements tangent to the identity with different orders then $G$ is not a  solvable group.
\end{Lemma}
\begin{proof}
Assume that there are $\hat f_{1}, \hat f_{2}\in G $ dicritic
diffeomorphisms of different orders, we say $p_{1} + 1$ and $p_{2} +
1$ respectively then by above lemma $\hat f_{3} = [\hat f_{1},\hat
f_{2}]$ is dicritic of order $p_{3} = p_{1} + p_{2} + 1 > p_{2} +
1$, similarly, we have that $\hat f_{4} = [\hat f_{3},\hat f_{2}]$
is dicritic of order $p_{4} = p_{3} + p_{2} > p_{3}$ and recurrently
$\hat f_{n} = [\hat f_{n - 1},\hat f_{n - 2}]$ is dicritic of order
$p_{n} = p_{n - 1} + p_{n - 2} > p_{n - 1}$, thus there is no
$n\in\mathbb{N}$ such that $ G  ^{(n)} = \{\Id\}$  and this
contradicts the fact $ G $ is solvable.
\end{proof}

\begin{proof}[Proof of Theorem~\ref{Theorem:D}]
Let $ G < \hat{\Diff}(\mathbb{C}^{n},0)$ be a subgroup of
diffeomorphisms tangent to the identity containing a dicritic
diffeomorphism $\hat f$ with order of tangency $k$. It is
immediate to verify that $(1) \Rightarrow (2)$. Let us now prove
$(2) \Rightarrow (3)$. Suppose that $\hat f(z) = z + f(z)z +
\cdots$. Suppose by contradiction that there is $\hat f^{(1)}\in G $
with order of tangency $p_{1} > k$ then obtain
\[
\hat f^{(2)} = [\hat f^{(1)},\hat f] = z + \hat f_{k_{2}}^{(2)} +
\cdots.
\]

We affirm that $\hat f_{k_{2}}^{(2)}\neq 0$ and thus $\hat f^{(2)}$
has order of tangent $k_{2} = k + k_{1} > k_{1} + 1$. In fact, as
the j-th coordinate of $\hat f_{k_{2}}^{(2)}$, is $(k_1 -
1)f.q_{k_1}^{(j)} - (\nabla f.Q_{k_1})z_j$, where $\hat f^{(1)} = z+
Q_{k_1} +\ldots$ and $Q_{k_1} = (q_{k_1}^{(1)},\ldots
,q_{k_1}^{(n)})$, $($ in consequence $k_{2} = k + k_{1})$, now if
$\hat f_{k_{2}}^{(2)} = 0$, then $(k_1 - 1)f.q_{k_1}^{(j)} = (\nabla
f.Q_{k_1})z_j$, for $j = 1,\ldots, n$. So following the same
argument of Lemma \ref{dicritico, qualq ----- dicrit} we have that
$Q_{k_1} = (g.z_1,\ldots,g.z_n)$, with $g$ homogeneous polinomial of
degree $p$, thus $Q_{k_1}$ has degree $k+1$, but this is impossible.
Repeating this process we can define:
\[
\hat f^{(n)} = [\hat f^{(n - 1)},\hat f] = z + \hat f_{k_{n}}^{(n)}
+ \cdots.
\]

Analogously $\hat f_{k_{n}}^{(n)}\neq 0$ and thus $\hat f^{(n)}$ we
has order of tangency $k_{n} = k + k_{n - 1} > k_{n - 1} + 1 > n$,
thus there is no  $n\in \mathbb{N}$, such that $\mathcal C ^{n}(G)
= \{Id\}$, what contradicts the fact that $G$ is nilpotent.
Therefore,
we have $G < \hat{\Diff}_{k+1}(\mathbb{C}^{n},0)$.\\
Now we prove $(3) \Rightarrow (1)$. From Lemma~\ref{comutador} we
have that for $\hat h, \hat g \in G $, $[\hat h,\hat g] = \{\Id\}$
or $[\hat h,\hat g]\in \hat{\Diff}_{\ell}(\mathbb{C}^{n},0)$,
$\ell\geq 2k+1$ thus $[\hat h,\hat g] = \{\Id\}$ and therefore $ G $
is abelian.
\end{proof}

\bibliographystyle{amsalpha}





\vglue.2in
\begin{tabular}{ll}
Mitchael Martelo  and   Bruno Scardua\\
Inst. Matematica - Universidade Federal do Rio de Janeiro\\
  Caixa Postal 68530 -  Rio de Janeiro-RJ\\
21.945-970 -  BRAZIL\\
michaelp204@gmail.com \, \\
scardua@im.ufrj.br
\end{tabular}

\end{document}